\documentclass[leqno,11pt]{amsart}
\usepackage{fancyhdr}

\usepackage[font=small,format=plain,labelfont=bf,up,textfont=it,up]{caption}
\usepackage{wrapfig}
\usepackage{graphicx}
\usepackage{rotating}
\usepackage{amsmath}
\usepackage{amssymb, latexsym}
\usepackage{amscd}
\usepackage{mathrsfs}
\usepackage{amsthm}
\usepackage{color}
\numberwithin{equation}{section}
\makeatletter
\def\tagform@#1{\maketag@@@{\textbf{(\ignorespaces#1\unskip\@@italiccorr)}}}
\makeatother

\begin{document}

\theoremstyle{plain}
\newtheorem{thm}{Theorem}[section]
\newtheorem{claim}[thm]{Claim}
\newtheorem{cor}[thm]{Corollary}
\newtheorem{prop}[thm]{Proposition}
\newtheorem{lemma}[thm]{Lemma}
\theoremstyle{definition}
\newtheorem{defn}[thm]{Definition}
\newtheorem{ex}[thm]{Exercise}
\newtheorem{prob}[thm]{Problem}
\newtheorem{remark}[thm]{Remark}

\title[ Homogeneous Subsets of a Lipschitz Graph and the Corona Theorem]{ \Large Homogeneous Subsets of a Lipschitz Graph and the Corona Theorem}

\author[Brady NewDelman]{Brady Max NewDelman\\  University of California, Los Angeles}
\address{University of California, Los Angeles}
\email{bnewdelm@math.ucla.edu}
\maketitle

\makeatletter
\def\section{\@ifstar\unnumberedsection\numberedsection}
\def\numberedsection{\@ifnextchar[
  \numberedsectionwithtwoarguments\numberedsectionwithoneargument}
\def\unnumberedsection{\@ifnextchar[
  \unnumberedsectionwithtwoarguments\unnumberedsectionwithoneargument}
\def\numberedsectionwithoneargument#1{\numberedsectionwithtwoarguments[#1]{#1}}
\def\unnumberedsectionwithoneargument#1{\unnumberedsectionwithtwoarguments[#1]{#1}}
\def\numberedsectionwithtwoarguments[#1]#2{%
  \ifhmode\par\fi
  \removelastskip
  \vskip 3ex\goodbreak
  \refstepcounter{section}%
  \begingroup
  \noindent\leavevmode\Large\bfseries\centering
  \S \thesection\ #2\par
  \endgroup
  \vskip 2ex\nobreak
  \addcontentsline{toc}{section}{%
    \protect\numberline{\thesection}%
    #1}%
  }
\def\unnumberedsectionwithtwoarguments[#1]#2{%
  \ifhmode\par\fi
  \removelastskip
  \vskip 3ex\goodbreak
  \begingroup
  \noindent\leavevmode\Large\bfseries\centering
  #2\par
  \endgroup
  \vskip 2ex\nobreak
  \addcontentsline{toc}{section}{%
    #1}%
  }

\newcommand{\HH}{\widetilde{\mathbb{H}}^+}

\bigskip \bigskip
\section{Introduction}

Let $\mathcal{A}: \mathbb{R} \rightarrow
\mathbb{R}$ be an $M$-Lipschitz continuous function.  Thus
$\mathcal{A}$ has a derivative almost everywhere such that $\Vert\mathcal{A}'\Vert_{L^{\infty}} = M$.  Let $\Gamma$
 be the Lipschitz graph parametrically defined
by $z(x) = x + i\mathcal{A}(x)$ in the extended complex plane, and
let $E_0$ be a closed set contained in $\Gamma$ with
\begin{equation}\label{density} 
 \Lambda\left(B(z,r) \cap E_0\right) >
 \epsilon_0 r \qquad\text{for all}~ z\in E_0, \quad \text{and all}~ r>0,
 \end{equation}
where $\Lambda$ is linear measure in the plane, $B(z,r)$ is
the open ball about $z$ of radius $r$, and $\epsilon_0>0$. The constant $\epsilon_0$ is called the \textit{Carleson lower density}. Any measurable subset of $\Gamma$ with a positive Carleson lower density is called
\textit{homogeneous} in $\Gamma$.\\

Set $\Omega = \mathbb{C}^* \backslash E_0$, and let $H^{\infty}(\Omega)$
 denote the space of bounded analytic functions on $\Omega$. In this paper, we prove:
\begin{thm}[The Corona Theorem] \label{thm}
Given $f_1, \ldots , f_n \in H^{\infty}(\Omega)$ and $\mu > 0$
with the property that $\mu \leq$ $\max \{|f_j(z)|: ~1\leq j \leq
n\} \leq 1$ for every $z \in \Omega$, there exist $g_1,
\ldots , g_n \in H^{\infty}(\Omega)$ such that $f_1g_1 + \cdots +
f_ng_n \equiv 1$ on $\Omega$.
\end{thm}
We will refer to the functions $\{f_j\}_{j=1}^n$ and
$\{g_j\}_{j=1}^n$ as the corona data and corona solutions respectively, and we will refer
 to $\mu$ and $n$ as the corona
constants.\\

There is an alternative way of viewing the theorem in the language
of uniform algebras. Let us denote by $\mathcal{M} =
\mathcal{M}(H^{\infty}(\Omega))$ the maximum ideal space of
$H^{\infty}(\Omega)$.  When $H^\infty(\Omega)$ separates the
points of $\Omega$, we can identify elements of $\Omega$ with
pointwise evaluation functionals in $\mathcal{M}$. Under this
identification, the theorem becomes equivalent to determining
whether $\Omega$ is dense in $\mathcal{M}$ in the Gelfand
topology. It is in this context where the theorem gets its name;
whereby, in the special case where $\Omega$ is the unit disk,
$\mathbb{D}$, we can think of $\mathbb{D}$ as being the sun, and
$\mathcal{M}
\backslash \overline{\mathbb{D}}$ as being the sun's corona.\\

Lennart Carleson (1962) proved the first corona theorem
for the case of the disk \cite{car-disc}.  His proof was subsequently simplified (using a $\overline{\partial}$ equation)
by H\"{o}rmander \cite{hormander},
 and later by a clever proof by Wolff (\cite{gamelin}, \cite{garnett}).
   The theorem was swiftly adapted to the case of finitely connected domains (Alling \cite{alling1}, \cite{alling2};  Stout \cite{stout1}, \cite{stout2}, \cite{stout3}; and others \cite{earle-marden}, \cite{forelli}, \cite{slod2}).  Each proof gave new insight into the structure of $H^{\infty}$.
The finitely connected domain proofs were fundamentally based upon
admixing localized corona solutions for
        overlapping simply connected components.
          One major drawback to this method was that the bounds of the corona solutions, $\Vert g_j\Vert_{\infty}$,
           were dependent on the number of boundary components.  This was an unfortunate hindrance as any planar domain can be exhausted by a sequence of finitely connected domains.
               Without a uniform bound on the corona solutions for the approximating domains, any method of taking
                normal limits was futile. In that direction, Gamelin \cite{gamelin2} observed that
                the corona theorem for all planar domains would be true if and only if
                there existed a uniform bound for finitely connected domains which is independent of the number of boundary components.  A proof to the corona problem for all planar domains remained a  mystery.\\

Further investigations into the corona problem revealed a connection between interpolating sequences,
 boundary thickness, and the Cauchy transform.  Along those lines, Carleson made another breakthrough by proving
  the corona theorem for domains with homogeneous boundary contained in the real line (homogeneous Denjoy domains).  The significance of his result was that these domains are infinitely connected.  Carleson lifted the corona data to the universal covering surface (where the corona solutions exist) and
      then mapped the solutions back to the original domain by an explicit projection operator invented by
       Forelli \cite{forelli}.  This concept was later simplified by Jones and Marshall \cite{jones-marshall}.  They determined
        that if the critical points of the Green's function for a domain form an interpolating sequence, then there exists a projection operator and the corona
         theorem is affirmative.  Moreover, they gave conditions necessary for determining when the critical points are
          indeed an interpolating sequence; one such condition can easily be proved when the boundary is
           homogeneous.  Following these results, the corona problem for all planar domains bounded by a
           homogeneous subset of a  graph seemed promising.\\

Peter Jones was the first to propose the idea of the corona problem for domains whose boundary lies in a Lipschitz graph \cite{jones3}.  He was motivated by the Denjoy conjecture, a consequence of Calder\'on's theorem on Cauchy integrals, which suggested that the space of bounded analytic functions was significantly abundant for these domains.  Thereby, one might be able to construct ``by hand'' the corona solutions.  As mentioned by Jones, the difficulty in the Lipschitz case was the lack of  symmetry. At that time, the deepest results for the corona theorem were in the Denjoy domains ($\Omega = \overline{\Omega}$) as in Carleson \cite{car-fatsets}, Jones and Marshall \cite{jones-marshall}, and Garnett and Jones \cite{garnettjones}.  These proofs made explicit use of the
  symmetry of the domains, either by confining the critical points to real intervals or by creating analytic functions by means of Schwarz reflection.  Nonetheless, Jones (unpublished) proved the corona theorem for domains bounded by a homogeneous subset of a Lipschitz graph.  He constructed by hand a projection operator akin to Forelli's.\\

For our proof, we work directly on the underlying space $\Omega$ without localizing the critical points of the Green's function, which can be cumbersome.   We divide $\Omega$ into two overlapping simply connected regions, $\widetilde{\Omega}^+$ and $\widetilde{\Omega}^-$.  On each region, we use Carleson's simply connected result to obtain regional corona solutions, $\{g_j^+\}_{j=1}^n$ and $\{g_j^-\}_{j=1}^n$.  Starting in $\widetilde{\Omega}^+$, we constructively solve a particular $\overline{\partial}$ equation to modify $\{g_j^+\}_{j=1}^n$  so that $\max_j |g_j^+(z)-g_j^-(z)|$ is reduced in the overlap of the regions.  After the modification, we do a similar procedure in $\widetilde{\Omega}^-$ to reduce the differences even more, then iterate the procedure.  The result of the iteration
gets us two uniformly bounded sequences of solutions on each
region. The $\overline{\partial}$ equation was constructed
specifically so that the normal limits of
the sequences agree on the overlap of the regions.\\

In proving the theorem, we assume that $E_0$ consists of a finite union of closed intervals
 in $\Gamma$, two of which are unbounded. This assumption is easily removed by a normal families argument
  provided that the number of intervals does
not control the bounds of the corona
solutions.\footnotetext[1]{The homogeneous condition combined with
the fact that $E_0$ is closed implies that $\epsilon_0 \leq 1/2$.
 The reason being that a complementary open interval $F \subset \Gamma \backslash E_0$ is not empty.
   As such, the double of F has a density less than $1/2$.} To be clear, when
we use the phrase, ``$J$ is an interval in $\Gamma$'' we mean
$p(J)$ is an interval in $\mathbb{R}$ for the projection $p:
\Gamma \rightarrow \mathbb{R}$ defined by $p(z(x)) = x$.
 It will also be convenient for us to consider
    the Lipschitz angle $\alpha = \tan^{-1}(M)$ for most of our calculations, instead of the slope $M$.\\

\noindent We mention here that there are two conditions equivalent to (1.1):

\begin{lemma}\label{lemmaequiv} When $\Gamma$ is an $M$-Lipschitz Graph and $E_0 \subset \Gamma$, the following three conditions are
equivalent:\\
\begin{enumerate}
\item[i)] \begin{quote} $E_0$ is homogeneous with a Carleson
    lower density $\epsilon_0$.\\  \end{quote}

\item[ii)] \begin{quote}There exists an $\epsilon_1 >0 $ such that
    $\left|p(E_0) \cap (x-r,x+r)\right| > \epsilon_1 r$
    \text{for all}~ $z=x+iy \in E_0$, \text{and all} $r>0$.\\  \end{quote}

\item[iii)]
 \begin{quote} There exists an $\epsilon_2 >0$ such that if we denote by $J_{z,r} = J_L \cup J_R$ the interval in
    $\Gamma$ containing $z$; $J_L$ is the subinterval having
    $z$ as a right
    endpoint, and $J_R$ is the subinterval
    having $z$ as a left endpoint with $\Lambda(J_L) =
    \Lambda(J_R) = r$, then $\Lambda(J_{z,r} \cap E_0) >
    \epsilon_2 r$, \text{for all} $z \in E_0$ \text{and all}
    $r>0$.\\  \end{quote}
\end{enumerate}

\noindent In addition, if either \emph{i)}, \emph{ii)}, or \emph{iii)} hold, then
\begin{enumerate}
\item[iv)] \begin{quote}There exists an $\epsilon_3>0 $ such that
    $\mathrm{cap} \left(B(z,r) \cap E_0\right)> \epsilon_3 r\\
    \text{for all}~ z \in E_0,~ \text{and all}~r>0 $.\\ \end{quote}
\end{enumerate}
\end{lemma}

The third item, iii), has the advantage that it applies to more
general curves, while iv) is even more general: it  says $E_0$ is
\textit{uniformly perfect} (see Pommerenke \cite{Pommerenke}). The
crux of the proof for Lemma 1.2. is based upon the relationship of
the projected length: \[\Lambda(J) \geq |p(J)| \geq
\cos(\alpha)\Lambda(J) \quad \text{for an interval $J \subset
\Gamma$.}\]\\

\noindent \textbf{Proof of Lemma 1.2:} \quad \noindent Let us first
assume that i) holds. Fix $z=x+iy \in \Gamma$ and $r>0$.
Since the projected mass of $B(z,r) \cap E_0$ lies inside of
$(x-r,x+r) \cap p(E_0)$ and
\[\left|p(B(z,r) \cap E_0)\right| \geq \cos(\alpha)  \Lambda\left(B(z,r) \cap E_0\right)
> \cos(\alpha)  \epsilon_0 r,\] \noindent  condition ii) holds with
$\epsilon_1 = \cos(\alpha)
\epsilon_0 $.\\

Now assume that ii) holds, and fix $z \in \Gamma$ and $r>0$.  By simple geometric considerations,
 we see that $B(z,r\cos(\alpha)) \cap \Gamma \subset J_{z,r}$.
 This implies
\begin{align*}
\Lambda( J_{z,r} \cap E_0) &\geq \Lambda\left(B(z,r\cos(\alpha)\right) \cap
E_0) \\
 &\geq |\left(x-r\cos(\alpha), x+r\cos(\alpha)\right) \cap p(E_0) | > \cos(\alpha)  \epsilon_1 r.
 \end{align*}
The last inequality is from ii). This implies condition iii) with $\epsilon_2 = \cos(\alpha)
\epsilon_1 $.\\

Showing that iii) implies  i) is simple as we can make the interval $J_{z, r}$ inside the ball $B(z,r)$. Then condition iii) implies $\Lambda(B(z,r)\cap E_0) >\epsilon_2 r$.  Thus $E_0$ is homogeneous with a Carleson lower density $\epsilon_2 $.\\

Lastly, from the proof of Theorem III.11 in Tsuji \cite{Tsuji}, we have the relationship for $E_0 \subset \Gamma$,
\[ \mathrm{cap}(B(z,r) \cap E_0) \geq \frac{\cos(\alpha) \Lambda(B(z,r) \cap E_0)}{2 e}.\]
 This tells us that i) implies iv) with $\epsilon_3 = \dfrac{\epsilon_0 \cos(\alpha)}{2 e}$. \qed \\
 
For the proof of Theorem 1.1, we make the additional assumption that the tangent
 to $\Gamma$ at a point $\zeta \in \Gamma$ is constant whenever $\zeta \in \Gamma \backslash E_0$. This comes
  without any loss of generality. Specifically, if we write  $\Gamma \backslash E_0 = \cup_k F_k ,$ then we define (see Figure 1.)
\[c_k = \tan^{-1}[\mathcal{A}'(x)], \qquad \emph{when}~\:z=x+iy \in F_k.\]
\begin{figure}[h]
\begin{center}
\includegraphics[width=1\textwidth , bb= 72 114 718 445]{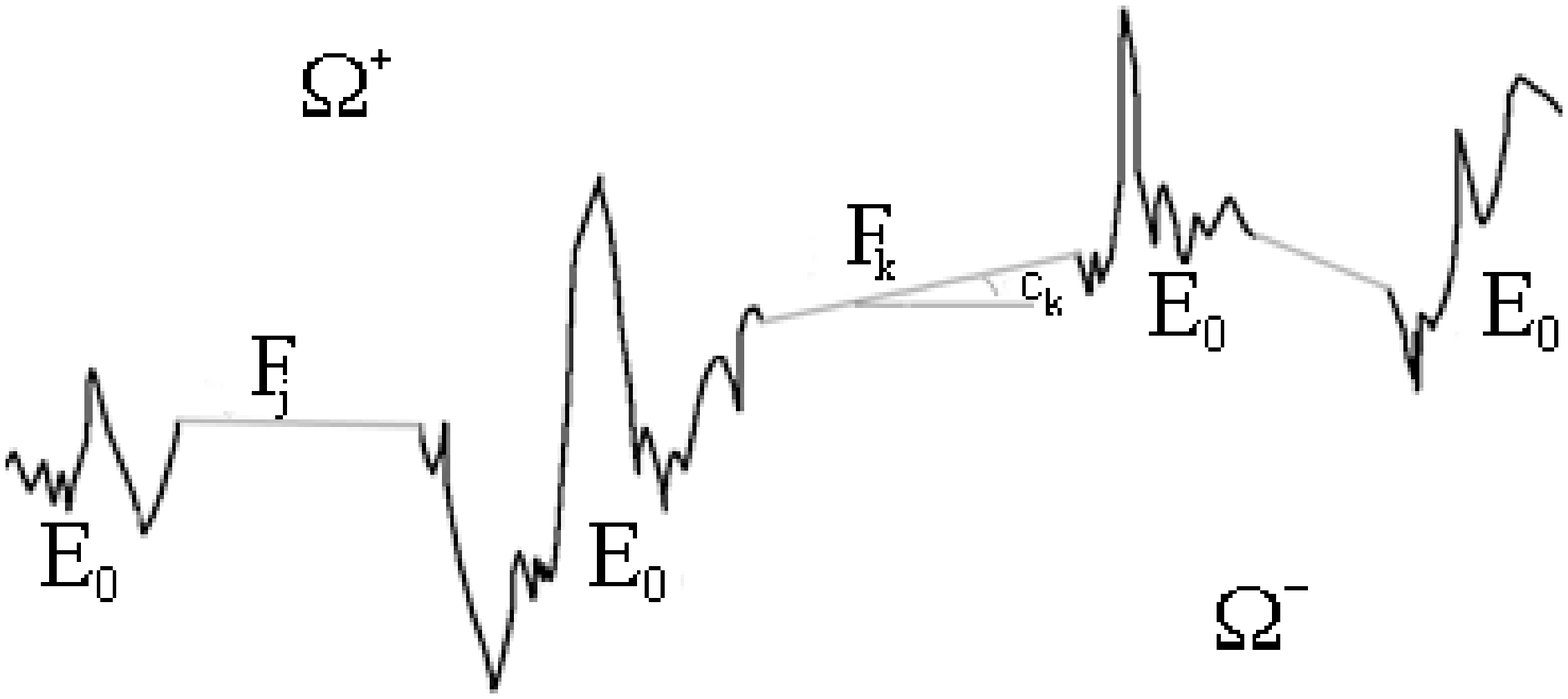}
\caption{The open interval $F_k$ makes an angle $c_k$ with the $x$-axis.}
\end{center}
\end{figure}

Let us now fix some notation that will be used throughout the
whole paper.  We define a tent region over an interval $J =
(z_1,z_2)$ with acute angle $\gamma$ by
\[ T_{(J,\gamma)}  = \left\{z: ~  0<\arg \left[ \frac{z-z_1}{z_2-z_1} \right]
 < \gamma \quad \emph{and}\quad 0< \arg \left[\frac{z_1-z_2}{z-z_2}\right]  < \gamma \right\}. \]

Notice that $\Gamma$ divides the plane into two simply connected
components, $\Omega^+$ and $\Omega^-$, where $\Omega^+$ lies above
$\Gamma$ and $\Omega^-$ lies below $\Gamma$.  With these
components, we fix two conformal maps and their inverses:
\begin{samepage}
\begin{align*}
 &\Phi^+(z) : \mathbb{H}^+ \rightarrow \Omega^{+}, \quad \Psi^+(z) =(\Phi^+(z))^{-1} : \Omega^+ \rightarrow \mathbb{H}^+ ,\\
  \intertext{\qquad \qquad \quad \emph{and}}
  & \Phi^-(z) : \mathbb{H}^- \rightarrow \Omega^{-}, \quad \Psi^-(z) =(\Phi^-(z))^{-1} : \Omega^- \rightarrow \mathbb{H}^- .\\
 \end{align*}
 \end{samepage}
We ask that $\Phi^+(\infty) = \infty$ and $\Phi^-(\infty) =
\infty$.  From Carath\'{e}odory's theorem we can extend our maps
to homeomorphisms so that $\Phi^+: \overline{\mathbb{H}^+}
\rightarrow  \overline{\Omega^+}$ and
 $\Phi^-: \overline{\mathbb{H}^-} \rightarrow  \overline{\Omega^-}$ respectively (see \cite{Kenig}, Theorem
 I.1.). \footnote[2]{By placing suitable minus signs, we may assume
 that $\mathfrak{Re}\{\Psi^{\pm}(z_1)\}<\mathfrak{Re}\{\Psi^{\pm}(z_2)\}$ whenever $z_1 ,z_2
 \in \Gamma$ and $\mathfrak{Re}\{z_1\} < \mathfrak{Re}\{z_2\}$.}
  We will be using two facts about about $\Phi^+$  and $\Phi^-$:

\begin{lemma}
The closed set $E^+ = \Psi^+(E_0)$  has a Carleson lower density $
\epsilon = \epsilon\left(\epsilon_0,\alpha\right) $ in
$\mathbb{R}$.  Likewise, $E^- = \Psi^-(E_0)$  has a Carleson lower
density $ \epsilon = \epsilon\left(\epsilon_0,\alpha\right)$ in
$\mathbb{R}.$
\end{lemma}

\begin{lemma}For any interval $F_j \subset \Gamma \backslash E_0$,
\begin{align*}
\Phi^+\left(T_{(\Psi^+(F_j),\gamma)}\right) &\subset T_{ (F_j,3 \gamma )} &\text{for~ $\gamma < \pi/12$}. \\
\text{Likewise},\quad \Phi^-\left(\overline{T_{(\Psi^-(F_j),\gamma)}}\right) &\subset \left(T_{ (F_j,3\gamma )}\right)^{*}
  &\text{for~ $\gamma < \pi/12$},
\end{align*}
where $*$ denotes reflection across $F_j$.\\
\end{lemma}

Lemma 1.3. tells us that homogeneity is preserved by the maps (although the Carleson lower densities may be different), while Lemma 1.4. tells
us that obtuse tents are mapped into obtuse tents.  It should be mentioned
 that $\pi/12$ is not crucial for Lemma 1.4. We made this choice since the acute angle gets tripled in the
 lemma and, throughout this paper, we will only consider tents that have an
 acute angle less than $\pi/4$.\\

\noindent \textbf{Proof of Lemma 1.3:}
  \quad We use a result of Kenig \cite{Kenig}:  if $\nu$ is the measure on $\mathbb{R}$ whose
density is $\left|(\Phi^+)^{\prime}(x)\right|$, then $\nu \in
\mathcal{A}_2$ on $\mathbb{R}$, where $\mathcal{A}_2$ is the class
of Muckenhoupt.  Now fix $r>0$, $x \in E^+$, $z =\Phi^+(x)$ and
have $I=(x-r,x+r)$.  Write $I = I_L \cup I_R$, where $I_L =
(x-r,x]$ and $I_R = [x,x+r)$ and denote $K=\Phi^+(I),$ $K_L =
\Phi^+(I_L)$, and $K_R = \Phi^+(I_R)$. Without loss of generality,
let us assume that $\Lambda(K_R) \leq \Lambda(K_L)$.  The
$\mathcal{A}_2$ relationship gives us a lower bound for
$\Lambda(K_R)$,
\begin{equation*}\label{1.2} \tag{1.2}
 \Lambda\left(K_R\right) \geq \frac{1}{C_2}\left(\frac{1}{2}\right)^2 \Lambda\left(K\right),
\end{equation*}
\noindent  where $C_2$ is the $\mathcal{A}_2$ constant.\\

Let $J_z$ be the interval inside $K$, as defined as in Lemma 1.2,
with $J_z = J_L \cup J_R$, where $J_R = K_R$ and $J_L$ is the
interval with right endpoint $z$ and length equal to
$\Lambda(K_R)$. From the proof of Lemma 1.2, we know that $\Lambda(E_0 \cap J_z) >
\epsilon_0 \cos(\alpha) \Lambda(K_R)$, and when we combine this inequality with
\eqref{1.2} we have
\begin{equation*}\label{1.3} \tag{1.3}
  \frac{\Lambda(E_0 \cap K)}{\Lambda(K)}   \geq  \frac{\Lambda(E_0 \cap J_z)}{\Lambda(K)} > \displaystyle \epsilon_0 \cos(\alpha)
   \frac{1}{C_2}\left(\frac{1}{2}\right)^2.
\end{equation*}

By a result of Muckenhoupt \cite{Muckenhoupt}, $\nu \in \mathcal{A}_2$ on $\mathbb{R}$  implies
 $\nu \in \mathcal{A}_{\infty}$ on $\mathbb{R}$. Hence, there exist constants $c_1>0$ and $c_2 >0$ independent of $E^+$ and $r$ such that,
\[ \left|\frac{E^+ \cap (x-r, x+r)}{(x-r,x+r)}\right| \geq c_1 \left(\frac{\Lambda\big(\Phi^+(E^+ \cap (x-r,x+r))\big)}
{\Lambda\big(\Phi^+((x-r,x+r))\big)}\right)^{\displaystyle c_2} =  c_1 \left(\frac{\Lambda(E_0 \cap K)}
{\Lambda(K)}\right)^{\displaystyle c_2}.\] Combining the above relationship with \eqref{1.3},
 we see that $E^+$ is homogeneous with a Carleson lower density depending only upon $\epsilon_0$ and $\alpha$. \qed \\

\noindent \textbf{Proof of Lemma 1.4:}  \quad The appearance of
the $*$ and the conjugation bar for the statement in the lower
half plane arise since the tents have an orientation to be above
the intervals.  It is not difficult to see that the two statements
remain alike upon modifying the arguments in the definition of the
tents, and we will only prove the result for the upper half
plane.\\

 Fix a tent domain $T_{ \left(I_j^+,\gamma\right)}$ over $
I_j^+ = \Psi^+(F_j) \in \mathbb{R}$ and write
$\log\left[(\Phi^+)'\right](z) = f_1(z) + i f_2(z)$ (take a
principle determination). Again from Kenig \cite{Kenig}, we have a
bounded argument for the derivative, that is $\left|f_2(z)\right|
\leq \alpha$ for all $z \in \mathbb{H}$. As such,  we can
represent $f_2(z)$ with a Poisson integral of the values coming
from its non-tangential limits on the real line:
\begin{align*}
f_2(z) - c_j &= \int_{\mathbb{R}} (f_2(t) - c_j)P_z(t)\;dt\\
&= \int_{I_j^+} (f_2(t) - c_j)P_z(t)\;dt + \int_{\mathbb{R} \backslash I_j^+} (f_2(t)-c_j)P_z(t)\;dt\\
&=  \int_{\mathbb{R} \backslash I_j^+} (f_2(t)-c_j)P_z(t)\;dt.
\end{align*}
\noindent The final equality holds since $f_2 = c_j$ over $I_j^+$. Taking absolute values of the above equality we get
 $|f_2(z) - c_j| \leq  2\alpha (1- \omega(z,I_j^+,\mathbb{H}^+)).$  Additionally, if
   $z \in T_{ (I_j^+,\gamma)}$ and $\gamma \leq \pi/12$, then \mbox{$\left|f_2(z) - c_j\right| \leq  4\alpha \gamma/\pi$~}
    by taking simple estimates for harmonic measure.  This lets us conclude that the values of the derivative lie in the
     cone domain:
\[c_j - \frac{4\alpha \gamma}{\pi} \leq \arg\left[\Phi^+(z) ^{\prime}\right] \leq c_j + \frac{4\alpha \gamma}{\pi} \qquad
 \emph{for ~$z\in T_{ \left(I_j^+,\gamma\right)}$.}\]

\noindent So that if we denote $I_j^+ = (x_1,x_2)$, then
\begin{align*}
\arg\left[\frac{\Phi^+(z) - \Phi^+(x_1)}{\Phi^+(x_2) - \Phi^+(x_1)}\right] &= \arg \left[ \frac{ \int_{[x_1,\, z]}~ \Phi'(w)\;dw}{\Phi^+(x_2) -\Phi^+(x_1)}\right] < \left(\gamma + \Big(c_j + \frac{4\alpha \gamma}{\pi}\Big)\right) - c_j < 3 \gamma,\\
\intertext{\emph{and}}
\arg\left[\frac{\Phi^+(x_1) - \Phi^+(x_2)}{\Phi^+(z) - \Phi^+(x_2)}\right] &= \arg \left[ \frac{\Phi^+(x_1) -\Phi^+(x_2)}
{ \int_{[x_2,\, z]}~ \Phi'(w)\;dw}\right] \\
&<\left(\pi + c_j\right) -\left( \left(\pi - \gamma\right) + \Big(c_j - \frac{4\gamma\alpha}{\pi}\Big) \right) < 3 \gamma.
\end{align*}

\noindent We conclude that $ \Phi^+(z)$ lies in $T_{ (F_j, 3 \gamma
)}$. \qed
\\
\\
\section{Four Crosscuts}

Recall $\Gamma \backslash E_0 = \cup_k F_k$, now let $\alpha_{M} =
(\pi/2 - \alpha)/4$, and let $\displaystyle D_{j}^+ = T_{ (F_j,
 \alpha_{M})}$  be the tent domain in $\Omega^+$ over $F_j$ with
acute angle $\alpha_{M}$, likewise define $D_{j}^- \subset
\Omega^-$. Merging the two tents together for all $j$, we make the
diamonds $D_j = D_{j}^+ \cup D_{j}^-$. The parameters for
$\alpha_M$ were chosen so that $\alpha_{M} < \pi/4$ and $D_j \cap
D_k = \emptyset$ for $j \neq k$ (see Figure 2).
\begin{figure}[h]
\begin{center}
\includegraphics[width=1\textwidth , bb= 72 144 718 460]{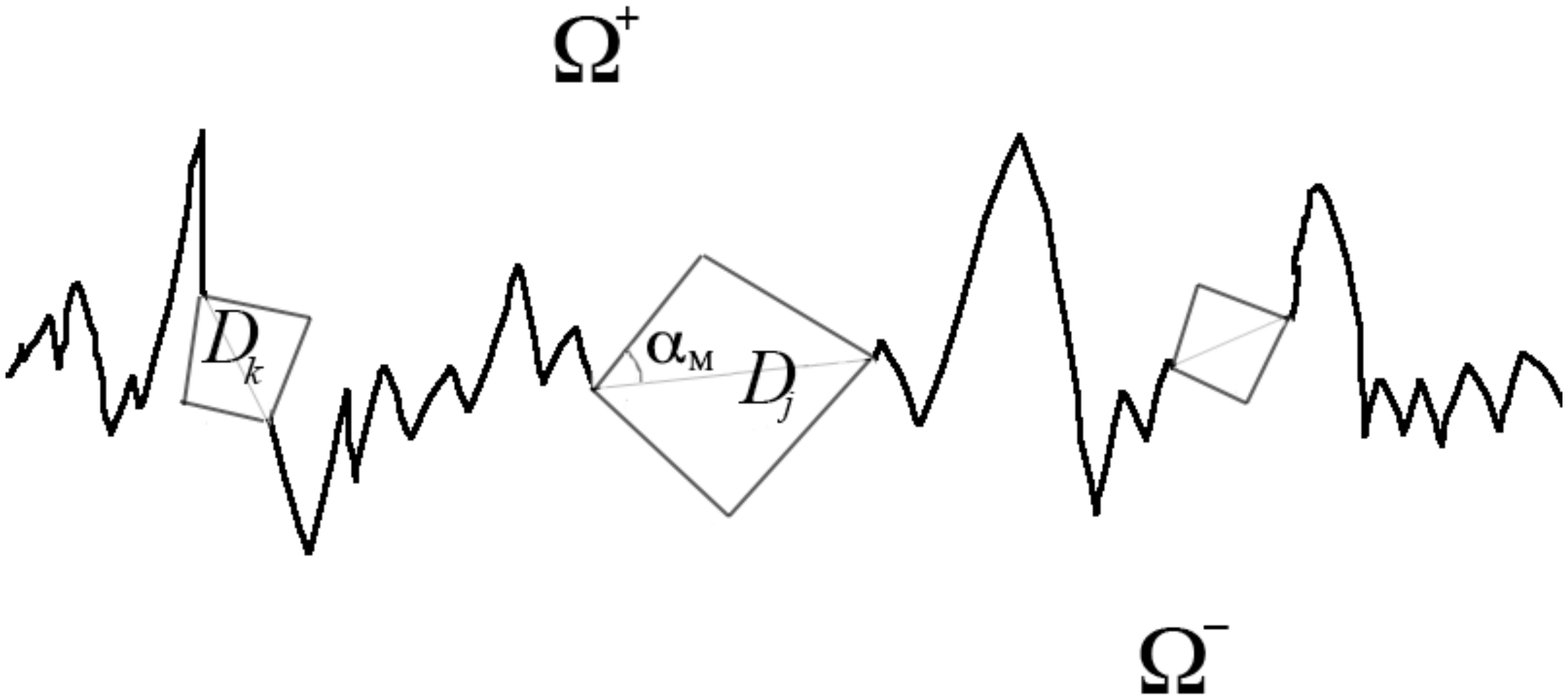}
\caption{\label{fig:a_figure} The diamond $\mathit{D_j}$ makes an acute angles $\alpha_M$ with the open interval $F_j$. The angle $\alpha_M$ is small enough to ensure that all diamonds are disjoint. }
\end{center}
\end{figure}

 In this section, we construct four families of crosscuts that encompass the open intervals $F_j$ and lie inside $D_j$.
   To do so, we will first need some elementary harmonic measure estimates.\\

In the upper half plane $\mathbb{H}^+$,
\[  \omega(z, E^+, \mathbb{H}^+) > \epsilon \qquad \emph{whenever} \quad z=x+iy, \quad \emph{and} \quad x \in E^+ .\] This is shown by decomposing the Poisson kernel for the upper half plane  into a sum of box kernels
 centered around $x \in E^{+}$ ($ P_z(t) = \sum_k a_k \: X_{A_k}(t)$). With this representation,
\begin{equation*} \omega(z, E^+, \mathbb{H}^+) =  \int_{E^+} P_z(t)\;dt =  \sum_k a_k|A_k \cap E^+|> \sum_k a_k  |A_k| \epsilon = \epsilon.
\end{equation*}
 If we write the complement of $E^+$ in the real line as  $\bigcup_j I_j^+ = \mathbb{R} \backslash E^+ = \bigcup_j \Psi^+(F_j)$, then $\omega(z, \cup I_j^+, \mathbb{H}^+) < 1-\epsilon$ when $\mathfrak{Re}\{z\} \in E^+$.  This tells us that we have the bounds $\omega(z, \cup I_j^+, \mathbb{H}^+)<1-\epsilon$ on the sides of the vertical half strip defined with base $I_j^+$ extending vertically in the upper half plane.  We can apply these bounds to get harmonic measure estimates on the boundaries of the diamonds:

\begin{lemma}
If  $z \in  \bigcup_j \partial T_{ (I_j^+,\gamma)}$ and $\gamma < \pi/4$, then $\omega(z, \cup I_j^+, \mathbb{H}^+) < 1- \dfrac{ \epsilon \gamma}{\pi}.$
\end{lemma}

\noindent \textbf{Proof of Lemma 2.1:}  \quad \noindent Fix $z \in \partial T_{ (I_k^+,\gamma)}$ for some $k$ and normalize $I_k^+$ into $(-\pi/2,\pi/2)$, and let us denote the half strip over $I_k^+$ by $S^+ = \left\{x+iy: ~ -\pi/2 < x < \pi/2  : ~ y>0\right\}$. From the preceding remarks,
\[ \omega(z, \cup I_k^+, \mathbb{H}^+) < (1-\epsilon) + \epsilon \: \omega(z, [-\pi/2,\pi/2], S^+).\]
If we write $ z = -\pi/2 + te^{i\gamma}$, then
\[ \arg\left[\sin(z) - \sin(-\pi/2)\right] = \arg\left[\int_{[-\pi/2,\, z]} \frac{d(\sin(w))}{dw}\; dw \right]
 = \arg \left[ \int_{[-\pi/2,\,z]} \cos(w)\; dw \right]\]
\[= \arg \left[ \int_0^t \cos(-\pi/2 + se^{i\gamma})e^{i\gamma} \; ds \right]
= \gamma + \arg \left[ \int_0^t \sin(se^{i\gamma})\; ds \right]  > \gamma. \] \\

\noindent By symmetry, this tells us that $\sin(z) \notin T_{([-1,1],\gamma)}$ when $z \in \partial T_{(I_k^+, \gamma)}$, so that
\[ \omega(z, [-\pi/2,\pi/2], S^+) = \omega(\sin(z), [-1,1], \mathbb{H}^+) < 1-\frac{\gamma}{\pi}.\]
Hence, \[ \omega(z, \cup I_k^+, \mathbb{H}^+) < (1-\epsilon) +
\epsilon \, (1- \frac{\gamma}{\pi}) =
 1 - \frac{\epsilon \gamma}{\pi}  \] \qed \\

We remark that Lemma 2.1. can easily be proved without conformal
maps but with a weaker bound on harmonic measure. This comes from
the observation that if we denote by $d= \textrm{dist}(z,E^+)$,
then $|E^+ \cap B(z,2d)| > \epsilon  d$. This implies for each $z
\in \bigcup_j \partial T_{(I_j^+,\gamma)}$ there exists a subset
of $E^+$ with linear measure proportionate to the distance of $z$
and the real axis.  The upper bounds for harmonic measure now
follow from estimating the Poisson kernel over these sets.
\\

With the estimates following from Lemma 2.1, we can now define our desired crosscuts.
  If we let $\beta_1 = 1 - \dfrac{\epsilon \alpha_{M}}{3 \pi}$ and
   $\beta_2 = 1 - \dfrac{1}{2} \dfrac{\epsilon \alpha_{M}}{3  \pi}$, then from Lemma 1.4. and Lemma 2.1,
\begin{align*}
&\gamma_{1}^{+} = \Phi^+\left(\left\{z:~ \omega(z, \cup I_j^+, \mathbb{H}^+)=\beta_1\right\}\right)  = \Phi^+\left(\delta_1^+\right)  \subset D_j,\\
\intertext{\qquad \qquad \emph{and}}
 &\gamma_{2}^{+} = \Phi^+\left(\left\{z:~ \omega(z, \cup I_j^+, \mathbb{H^+})=\beta_2\right\}\right) = \Phi^+\left(\delta_2^+\right)  \subset D_j.
\end{align*}
Similarly, we define the $\gamma_1^-$, $\gamma_2^-$, $\delta_1^-$, and $\delta_2^-$ for the lower half planes. These will be our collection of crosscuts.\\

Recall, a \textit{Carleson contour} in the upper half plane is a
countable union $\mathcal{C}$ of rectifiable arcs in
$\mathbb{H^+}$ such that for every interval $I \subset \mathbb{R}$,
\[ \Lambda \left(\mathcal{C} \cap (I \times (0,|I|)\right) \leq
C(\mathcal{C})\: |I|.\] This implies arc length on $\mathcal{C}$ is a Carleson measure with constant $C(\mathcal{C}).$
\begin{lemma} The crosscuts
 $ \delta_1^+ = \left\{z:~ \omega(z, \cup I_j^+, \mathbb{H}^+)=\beta_1\right\} $ form a Carleson contour in $\mathbb{H}^+$.
  Likewise, $\delta_1^- = \left\{z:~ \omega(z, \cup I_j^-, \mathbb{H}^-)=\beta_1\right\}$ form a Carleson contour in $\mathbb{H}^-$.
\end{lemma}

\noindent \textbf{Proof of Lemma 2.2:}  \quad First we recall that
$\delta_1^+ = \{z: \omega(z,E^+,\mathbb{H}^+) = 1-\beta_1\}$ lies
under the tents $\bigcup_{k}
T_{\left(I_k,\frac{\alpha_M}{3}\right)}$.  Next, if we write $E^+
= \bigcup_{j} [a_j, b_j]$, then \footnotetext[3]{We use the
convention $\arg[\infty,z] = 0$ and $\arg[-\infty,z] = \pi$
respectively when $[a_j, b_j] = [a_j, \infty]$ and $[a_j, b_j] =
[-\infty, b_j]$. }
\begin{align*}
\omega(z,E^+, \mathbb{H}^+) &= \frac{1}{\pi} \sum_j \arg \left[\frac{b_j -z}{a_j -z} \right],\\
\intertext{and by taking a derivative,}
\omega_y(z) + i\omega_x(z) &= \frac{1}{\pi} \sum_j \left(
\frac{1}{z-b_j} - \frac{1}{z-a_j} \right).\\
\end{align*}
Separating the real and imaginary parts gives us the ratio\\

\begin{equation*}\tag{2.1} \label{2.1}
 \displaystyle \dfrac{\omega_x(z)}{\omega_y(z)} = \dfrac{\displaystyle \sum \frac{(b_j-a_j)~\mathfrak{Im}\{\overline{(z-b_j)(z-a_j)}\}}{|z-a_j|^2 |z-b_j|^2}       }{\displaystyle \sum \frac{(b_j-a_j) ~\mathfrak{Re}\{\overline{(z-b_j)(z-a_j)}\}}{|z-a_j|^2 |z-b_j|^2} }.
\end{equation*}
\\

Suppose $z \in \bigcup_k T_{\left(I_k,\frac{\alpha_M}{3}\right)}$, then $|\arg
\left[(z-b_j)(z-a_j)\right]| < 2\alpha_M/3$ for all $j$; and since
$\frac{2\alpha_M}{3} < \frac{\pi}{4}$, this makes
\begin{equation*}\tag{2.2}\label{2.2}
 \frac{\left|\mathfrak{Im}\{(z-b_j)(z-a_j)\}\right|}{ \mathfrak{Re}\{(z-b_j)(z-a_j)\}}
    \leq \tan\left(\frac{2\alpha_M}{3}\right) \quad \emph{for
    all $j$}.
\end{equation*}

\noindent It is also clear that $\omega_y(z)>0$ for all $z$ in the
tents, so that by comparing the like terms in the sums of (2.1)
with the ratio in (2.2),
\[ \left| \frac{\omega_x(z)}{\omega_y(z)} \right| \leq \tan\left(\frac{2\alpha_M}{3}\right). \]
As the curve $\delta_1^+$ is a level set, the gradient of $\omega$ at any point is perpendicular to the tangent of the curve.  With the above ratio, we conclude that the tangent to the level curves is bounded in argument by $\frac{2\alpha_M}{3}$.
  This means that $\delta_1^+$ is a Carleson curve with a constant of $\sec\left(\frac{2\alpha_M}{3}\right)$. \qed
\\
\\

\section{The Regions $\mathcal{D^+}$ and $\mathcal{D^-}$}

 Now that we have the cross cuts $\{\gamma_{j}^{\pm} \}_{j=1,2}$,
we may define
 the following extended domains: let $\widetilde{\Omega}^+$ be the simply connected domain
  containing $\Omega^+$ that is bounded by the closed intervals of $E_0$ and the bottom crosscuts $\gamma_1^-$.
   As $\Psi^+(z)$ has a constant argument on each $F_j$, and the crosscuts of $\gamma_1^-$ lie in  disjoint diamonds, $\Psi^+(z)$
    can be extended (by reflecting across each $F_j$) to a
    map $\widetilde{\Psi}^+: \widetilde{\Omega}^+ \rightarrow \widetilde{\mathbb{H}}^+$,
     where $\widetilde{\mathbb{H}}^+$ is the domain containing $\mathbb{H}^+$
      that is bounded by $E^+$ and
      $\widetilde{\Psi}^+(\gamma_1^-).$  \footnote[4]{We are not identifying the crosscuts $\widetilde{\Psi}^+(\gamma_1^-)$ with
       the crosscuts $\delta_1^-$.}   \\

 \subsection*{Interpolating Functions}
\indent \indent A sequence $\{z_m\}_{m=1}^{\infty} \subset \mathbb{H}^+$ is called an \textit{interpolating sequence} for $H^{\infty}(\mathbb{H}^+)$ if, whenever $|w_m| \leq 1$, there exists a function $f\in H^{\infty}(\mathbb{H}^+)$ such that
\[f(z_m) = w_m, \quad m=1,2,\ldots\]
 When $\{z_m\}_{m=1}^{\infty}$ is an interpolating sequence, we call the finite bound
\[ \mathcal{N}\left(\{z_m\},\mathbb{H}^+\right) = \sup_{|w_j| \leq 1} \inf \left\{\|f\|: ~ f \in H^{\infty}(\mathbb{H}^+) \quad \emph{and} \quad f(z_m)=w_m, ~~ m=1,2,\ldots \right\}\]
the \textit{constant of interpolation}. By a theorem of Carleson, $\{z_m\}_{m=1}^{\infty}$ is an interpolating sequence
 if and only if
\[ \delta_{\mathbb{H}^+} \left(\{z_m\} \right) = \inf_n  \displaystyle \prod_{~k, \, k \neq n}
 \left| \frac{z_n - z_k}{z_n - \overline{z_k}} \right| >0;\]
furthermore, we have the relationship $1/\delta_{\mathbb{H}^+} \leq \mathcal{N} \leq ~ (1-\log \delta_{\mathbb{H}^+})c/\delta_{\mathbb{H}^+}$, in which $c$ is some absolute constant.  For a nice discussion on interpolating sequences and a proof of Carleson's interpolation theorem see \cite{garnett}, \S VII.\\

Fix $A = \dfrac{1-\beta_1}{1+3\beta_1}$, and let $\{z_m\}_{m=1}^{\infty} \subset \mathbb{H}^+$ be a sequence embedded in $\delta_1^+$ satisfying
 \begin{equation*}
 |z_n - z_m| \geq Ay_m, \quad n \neq m.
 \end{equation*}
 Since the sequence lies in a Carleson contour (by Lemma 2.2) and it is hyperbolically separated,
   we know that $\delta_{\mathbb{H}^+}(\{z_m\}) = C(A, \epsilon, \alpha_M)>0$
    (see \cite{garnett}, \textbf{\S}\textrm{VII}). Carleson's interpolation theorem then implies
     that $\{z_m\}_{m=1}^{\infty}$ is an interpolating sequence for $\mathbb{H}^+$. It is
     also the case that $\{z_m\}_{m=1}^{\infty}$ is an interpolating sequence for the extended domain $\widetilde{\mathbb{H}}^+$.
This follows from Garnett and Jones \cite{garnettjones} (Theorem IV.1), and applies in our case since
$\widetilde{\mathbb{H}}^+$ is a subset of a Denjoy domain.
Alternatively, $\{z_m\}_{m=1}^{\infty}$ can be shown to be an interpolating
sequence for $\widetilde{\mathbb{H}}^+$ by a result of
Gonz\'{a}lez and Nicolau \cite{gonzaleznicolau}. From their
result, it suffices to have $\delta_{\mathbb{H}^+}>0$ for the
image of $\{\Phi^+(z_m)\}$ under the canonical quasi-conformal map
that takes the domain $(\cup_j D_j) \cup \Omega^+ $ to the upper half
plane. Since the quasi-conformal map is explicit, it is easy to
verify; we omit the details. In any case, there exist
interpolating functions for the domain
         $\widetilde{\mathbb{H}}^+$ with a constant of interpolation $\mathcal{N} =  \mathcal{N}(A, \alpha_M, \epsilon)$.\\

Working again in the upper half plane, let $ B =\min\left\{A,\:  \dfrac{1}{6\mathcal{N}^2}\right\}$
 and denote the region \[
\mathcal{D}^+ = \left\{z \in \mathbb{H}^+: ~
 \omega(z, \cup I_j^+, \mathbb{H}^+) > \beta_1,
~ \: d(z) < B \right\},\] with $d(z) = y^{-1} \inf_{\zeta \in
\delta_1^+}|z-\zeta|$.  We chose $A$ so that by Harnack's
inequality $\mathcal{D}^+$ lies above $\delta_2^+$, \mbox{that is}
 \[ \omega(z,\cup I_j^+, \mathbb{H}^+) < \beta_2,
\quad \emph{for all}~ z \in
\mathcal{D}^+ \qquad (\text{see Figure 3.}).\]
\begin{figure}[h]
\begin{center}
\includegraphics[width=1\textwidth , bb= 72 164 718 445]{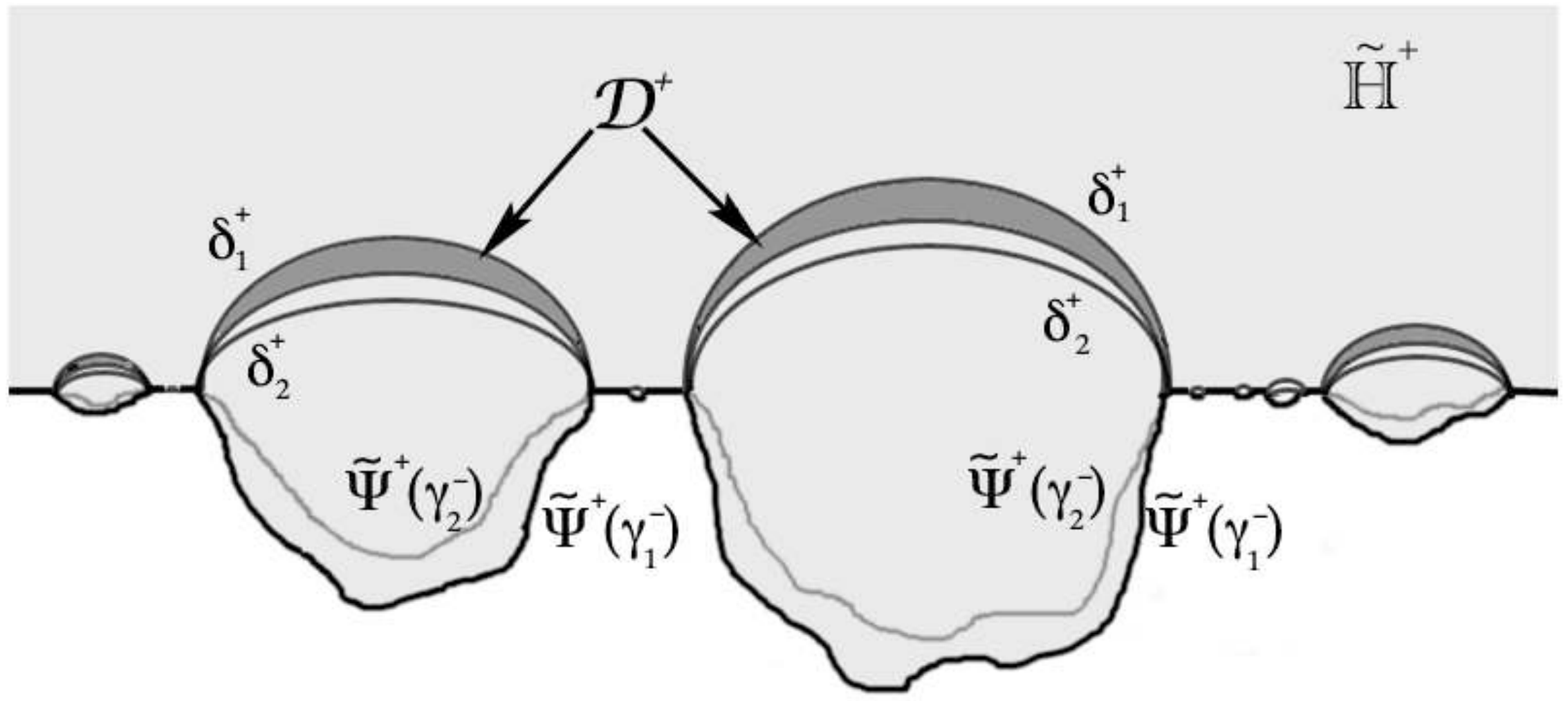}
\caption{The region $\mathcal{D}^+$ lies between the curves $\delta_1^+$ and $\delta_2^+$ and separates the two extended half planes.}
\end{center}
\end{figure}

\noindent Now fix a sequence $\{z_m^+\} \subset \delta_1^+$ satisfying
  \begin{align}
 |z_n^+ - z_m^+| &\geq By_m^+, \qquad n \neq m, \label{3.1} \\
\inf_m \frac{|z-z_m^+|}{y_m^+} &\leq 3B, \qquad \emph{for all} ~z \in \mathcal{D}^+.
 \end{align}
The existence of such a sequence follows by taking a maximal sequence
  satisfying \eqref{3.1}.  We now follow a standard argument as originated in Garnett and Jones \cite{garnettjones} (Lemma II.2) and as used in Handy \cite{Handy} (Lemma III.2) to obtain a specific set of interpolation functions:

\bigskip

\begin{samepage}
\begin{lemma}\label{interpolating sequences}
There exists functions $\{h_m^+\}_{m=1}^{\infty} \subset H^{\infty}(\widetilde{\mathbb{H}}^+)$ such that
\begin{align}  h_m^+(z_m^+) &= 1, \hspace{.5in} \label{3.3} \\
 \left\Vert h_m^+ \right\Vert_{H^{\infty}} &\leq  \mathcal{N}^2, \hspace{.5in} \label{3.4} \\
 \textit{and} \qquad \qquad \qquad \quad& \notag \\
\sum_m \left|h_m^+(z)\right| &\leq\mathcal{K}(A, \epsilon, \alpha_M) \quad z\in \widetilde{\mathbb{H}}^+.  \hspace{.5in} \label{3.5}
\end{align}
\end{lemma}
\end{samepage}

\noindent \textbf{Proof of Lemma 3.1:} \quad By a stopping time argument used to group $\{z_m^+\}$ into generations
 (see \cite{garnett}, pg. 416), we may split $\{z_m^+\}$ into a finite union of disjoint subsequences
  $S_k$, $1 \leq k \leq 2^p$, so that
\[   \inf \left\{  \frac{\left|z_j^+ - z_l^+\right|}{y_l^+}: ~z_j^+,\, z_l^+ \in S_k, \quad j\neq l \right\} \geq A. \qquad \emph{for all $k$.}\]

\noindent Since the points of $S_k$ are hyperbolically separated by $A$, our earlier discussion implies that each $S_k$ has a constant of interpolation less than $\mathcal{N}$.\\

Let us restrict our attention to a fixed subsequence $S_k$.  If we assume that $S_k = \{z_1, z_2, \ldots, z_{n_0}\}$ is
 finite, then there exists $f_j \in H^{\infty}(\widetilde{\mathbb{H}}^+)$ such
  that $\Vert f_j\Vert_{H^{\infty}} \leq \mathcal{N}$ and $f_j(z_m) = \omega^{mj}$, where
   $\omega= e^{2\pi i/n_0}$.  Moreover, if we define
\begin{align*}
 h_m^+(z) &= \left(\frac{1}{n_0} \sum_{j=1}^{n_0} \omega^{-mj}f_j(z) \right)^2 ,\\
\intertext{then $h_m^+(z_j) = \delta_{m,j}$ and}
\sum_{m=1}^{n_0} |h_m^+(z)| &= n_0^{-2}\sum_{m=1}^{n_0} \sum_{j,l} \omega^{-mj} \omega^{ml}f_j(z)\overline{f_l}(z) \\
&= n_0^{-2} \sum_{j=1}^{n_0} n_0|f_j(z)|^2 \leq \mathcal{N}^2.
\end{align*}
Therefore, by exhausting each $S_k$ and taking the normal limits, we have \eqref{3.3}, \eqref{3.4}, and \eqref{3.5} with $\mathcal{K} = \mathcal{N}^2 2^p.$ \qed \\

The technique used above of averaging interpolating functions is due to Varopoulos \cite{varopo}.  We made our choice of $\mathcal{N}^2B\leq 1/6$ specifically so that if we write $\mathcal{D}^+$ as the disjoint union of sets
$\mathcal{D}_n^+ \subset \{z: ~ |z-z_n^+| \leq
3By_n^+\}$, then with \eqref{3.4} and Schwarz Lemma

\begin{equation}\label{3.6}
|h_n^+(z)| > 1/2 \qquad \emph{whenever} ~ z \in \mathcal{D}_n^+.
\end{equation}

Since throughout this chapter we could change all plus signs to minus signs, we could likewise define our friends: $\widetilde{\Omega}^-$, $\widetilde{\Psi}^-$, $\widetilde{\mathbb{H}}^-$,  $\mathcal{D}^-$, $\{z_m^-\}_{m=1}^{\infty}$, and $\{h_m^-\}_{m=1}^{\infty}$.
\\
\\

\section{Iterative Blending of Corona Solutions}
We now begin the process of ``sewing'' together corona solutions
from the simply connected domains $\widetilde{\Omega}^+$ and
$\widetilde{\Omega}^-$. Let $\{g_j^0\}_{j=1}^n$ be an arbitrary
corona solution set for $\widetilde{\Omega}^+$  and let
$\{g_j^1\}_{j=1}^n$ be an arbitrary corona solution set for
$\widetilde{\Omega}^-$.  These solution sets exist from Carleson's
simply connected corona theorem; furthermore, there is a uniform
bound for the sets:
 \[ \left\Vert g_j^0 \right\Vert_{H^{\infty}(\widetilde{\Omega}^+)} \leq N  \quad \emph{and} \quad \left\Vert g_j^1\right\Vert_{H^{\infty}(\widetilde{\Omega}^-)}  \leq N, \quad  j=1, \ldots ,n.\]
 The bound, $N$, depends only on
 the corona constants: $N = N(\mu,\delta, n)$  (\cite{garnett} \textbf{\S}IIX).  In
  this chapter, we are going to create a special collection of solutions
   $\{g_j^k\}_{j=1}^n \subset H^{\infty}(\widetilde{\Omega}^+)$ when $k$ is even,
and $\{g_j^k\}_{j=1}^n \subset
H^{\infty}(\widetilde{\Omega}^-)$ when $k$ is odd.\\

\subsection*{The First Stitchings}
\indent \indent The first sewing of the corona solutions will be
across the region $\mathcal{D}^+$ in $\widetilde{\mathbb{H}}^+$.
Let us denote by $\omega(z) = \omega(z, \cup I_j^+, \mathbb{H}^+)$
as the harmonic measure for the open intervals $\bigcup_j I_j^+$
in the upper half plane. Although $\omega(z)$ is not defined on
the extended domain, we use the convention $\omega(\bar{z}) = 2 -
\omega(z, \cup I_j^+, \mathbb{H}^+)$ to extend $\omega$ to
$\widetilde{\mathbb{H}}^+$.  (When working in
$\widetilde{\mathbb{H}}^-$, we will also be denoting with
$\omega(z)$ and it should be clear from context.)\\

We can think of both $\{f_j(z)\}_{j=1}^n$ and
$\{g_j^0(z)\}_{j=1}^n$ as being defined in
$\widetilde{\mathbb{H}}^+$ under the map $\widetilde{\Phi}^+(z)$,
and we will not change our notation. On the other hand, when we
write $\{g_j^1(z)\}_{j=1}^n$ we must remember that these functions
are only defined in $\widetilde{\mathbb{H}}^+$ between the curves
$\{\widetilde{\Psi}^+(\gamma_1^-)\}$ and $\{\delta_1^+\}$ (see
Figure 4).  Because of a calculation advantage, we have chosen not
to sew on the $\widetilde{\Omega}^+$ side where the corona data
and solutions are originally defined, but instead work in the
extended half planes where we have defined $\mathcal{D}^+,
\{z_m^+\},$ and $\{h_m^+(z)\}$.\begin{figure}[h]
\begin{center}
\includegraphics[width=1\textwidth , bb= 72 174 718 445]{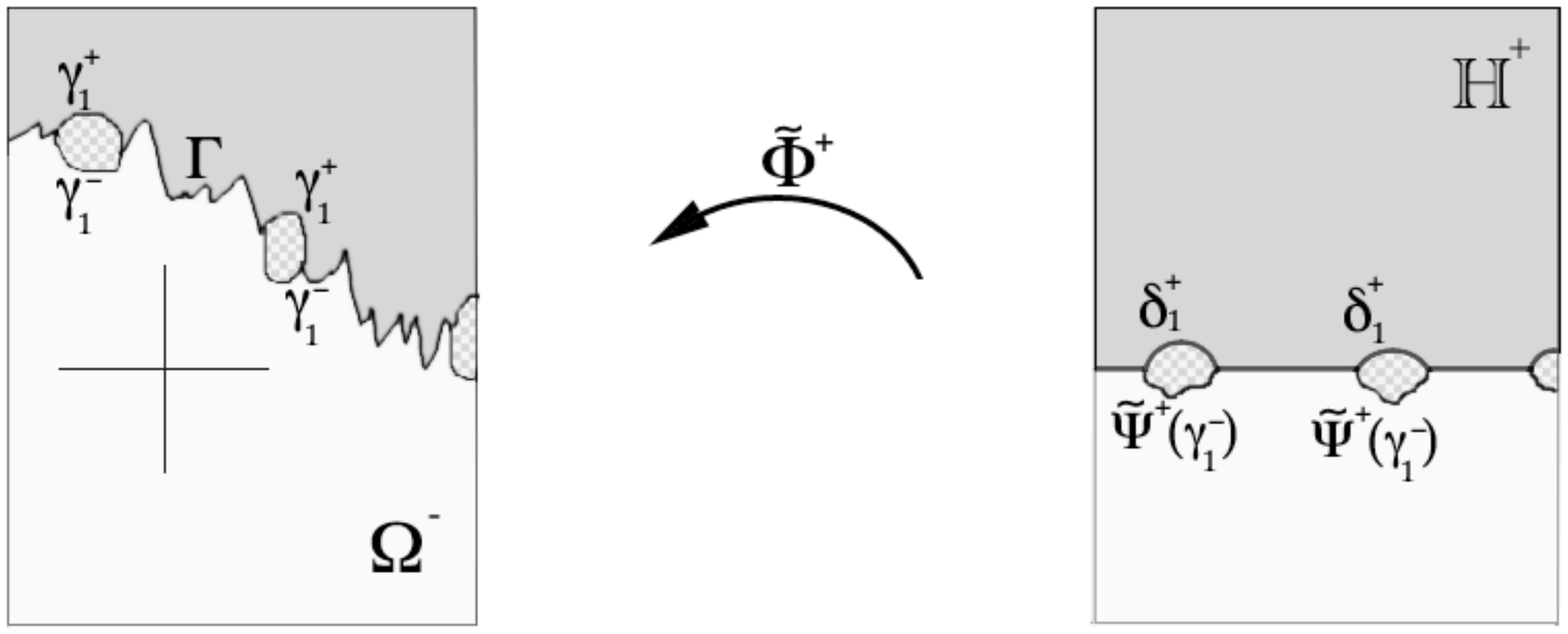}
\caption{The checkered regions lie in the intersection of the domains for $\{g_j^0\}_{j=1}^n$ and $\{g_j^1\}_{j=1}^n.$}
\end{center}
\end{figure}

Let $\{\varphi^+(z), 1- \varphi^+(z) \}$ be a smooth partition of
unity
 for $\widetilde{\mathbb{H}}^+$ across $\mathcal{D}^+$, with
$\varphi^+(z) =1$ when $\omega(z) \leq \beta_1$.  By standard
arguments, $|y_n^+|| \nabla \varphi^+(z)| \leq C B^{-1}$ for all
$z \in \mathcal{D}_n^+$, where $C$ is a positive constant.   Using
our partition, let us piece together the two families of corona
solutions: \[G_j^2 = g_j^0 \varphi^+ + g_j^1(1-\varphi^+) \qquad
j=1,\ldots ,n.\] These smooth functions are a well defined
solution set to the corona equation \[\sum f_j(z)G_j^2(z)=1\] for
the region $\widetilde{\mathbb{H}}^+$, but they are not
necessarily analytic. Therefore, we consider a technique of
H\"ormander \cite{hormander} (and as used in \cite{garnettjones}
and \cite{Handy}).  We seek to find functions \mbox{$\{a^2_{j,k}\}
\subset L^{\infty}(\widetilde{\mathbb{H}}^+)$} that solve (in the
sense of distributions) the $\overline{\partial}$ equation
\[ \overline{\partial} a_{j,k}^2 = G_j^2 ~ \overline{\partial} G_k^2.\]
Indeed, such functions provide the necessary cancelation to make
the collection
\begin{equation}\label{4.1}
g_j^2 = G_j^2 + \sum_{k=1}^n \left(a_{j,k}^2 - a_{k,j}^2\right)f_k \qquad  j=1,\ldots, n
\end{equation} a solution set while simultaneously solving the equation $\overline{\partial} g_j^2 =0$ in the sense of distributions.  Then upon modifying each $a_{j,k}^2$ on a set of measure zero, Weyl's lemma will allow us to conclude that the collection $\{g_j^2\}_{j=1}^n$ is a bona fide corona solution set in $\widetilde{\mathbb{H}}^+$.\\

For our construction, we require not only that the functions
$\{a_{j,k}^2\}$ are bounded, but also have an additional
convergence factor.  Fix $1> b_1 >0$ (to be determined later) and
denote by $\tilde{\omega}(z)$ as the harmonic conjugate for
$\omega(z)$.  Consider the equation
 \begin{align*}
  a_{j,k}^2(z) &= \frac{1}{\pi} \sum_l
    \iint_{\mathcal{D}_l^+} \left(b_1^{ \omega(z) -
    \omega(\zeta)
  + i \left(\widetilde{\omega}(z) -
   \widetilde{\omega}(\zeta)\right)}\right)~ \frac{G_j^2(\zeta)
   \overline{\partial} G_k^2(\zeta)}{\zeta -z}
   \frac{h_l^+(z)}{h_l^+(\zeta)}\; d\zeta d\bar{\zeta}.\\
   \intertext{Formally $\overline{\partial} a_{j,k}^2 = G_j^2 ~ \overline{\partial} G_k^2$, so we need to check the convergence of the sum}
\left |a_{j,k}^2(z)\right| &\leq \frac{1}{\pi}\sum_l \iint_{\mathcal{D}_l^+}
 \left(b_1^{ (\omega(z) - \omega(\zeta) }\right)~ \frac{|G_j^2(\zeta) \overline{\partial}
  G_k^2(\zeta)|}{|\zeta -z|} \frac{|h_l^+(z)|}{|h_l^+(\zeta)|}\; d\zeta d\bar{\zeta}.  \\
\intertext{Using \eqref{3.6} and recalling $\omega(\zeta)<\beta_2$ when $z\in \mathcal{D}^+$,}
 \left |a_{j,k}^2(z)\right|&\leq \frac{2}{\pi} \sum_l |h_l^+(z)|\iint_{\mathcal{D}_l^+} \left(b_1^{ \omega(z) - \omega(\zeta) }\right)~
 \frac{|G_j^2(\zeta) \overline{\partial} G_k^2(\zeta)|}{|\zeta -z|} \; d\zeta d\bar{\zeta}  \\
&\leq \frac{2}{\pi} \sum_l |h_l^+(z)| \left(b_1^{ (\omega(z) - \beta_2)}\right) \iint_{\mathcal{D}_l^+}
   \frac{|G_j^2(\zeta) \overline{\partial} G_k^2(\zeta)|}{|\zeta -z|} \; d\zeta d\bar{\zeta}  \\
 &\leq \frac{2}{\pi} \sum_l |h_l^+(z)| \left(b_1^{(\omega(z) - \beta_2)}\right) \iint_{\mathcal{D}_l^+}
  \frac{|G_j^2(\zeta)| |g_k^1(\zeta) - g_k^0(\zeta)| |\nabla\varphi^+(\zeta)|}{|\zeta -z|}\;  d\zeta d\bar{\zeta}.
\end{align*}
Before we show the above is a convergent sum, we would like to
identify some key numbers that will appear in the iterative
process.  Using the notation $\left\Vert \cdot
\right\Vert_{\mathcal{D}^+}$ and $\left\Vert \cdot
\right\Vert_{\mathcal{D}^-}$ for the supremum of the modulus in
the region $\mathcal{D}^+$ and $\mathcal{D}^-$ respectively,  let
us label

\begin{align*}
x_m &=
\begin{cases} \displaystyle \max _k \left\Vert g_k^{m} - g_k^{m-1} \right\Vert_{\mathcal{D}^+}    & \emph{when $m$ is odd,}
\\
\displaystyle \max_k   \left\Vert g_k^{m} - g_k^{m-1} \right\Vert_{\mathcal{D}^-}       &\emph{when $m$ is even,} \end{cases} \nonumber \\
y_m &=
\begin{cases}\displaystyle \max_k  \left\Vert G_k^{m+1} \right\Vert_{\mathcal{D}^+}     & \emph{\quad when $m$ is odd,}
\\
\displaystyle \max_k   \left\Vert  G_k^{m+1}  \right\Vert_{\mathcal{D}^-}       &\emph{\quad when $m$ is even.}
\end{cases} \nonumber
\end{align*}

\noindent So that in our context,

\begin{align*}
 \left|a_{j,k}^2(z)\right| &\leq \frac{2}{\pi}  \sum_l \left|h_l^+(z)\right| b_1^{(\omega(z) - \beta_2) } x_1 y_1 \iint_{\mathcal{D}_l^+}  \frac{|\nabla\varphi^+(\zeta)|}{|\zeta -z|}\;  d\zeta d\bar{\zeta}   \\
   &\leq \left(12 C\right) \sum_l |h_l^+(z)| b_1^{ (\omega(z) - \beta_2) } x_1 y_1, \\
\intertext{and by using \eqref{3.5} we reduce the inequality to}
 \left|a_{j,k}^2(z)\right| &\leq  (12C)  ~ \mathcal{K} ~ b_1^{\left(\omega(z) - \beta_2\right) }~ x_1 y_1.
\end{align*}

We conclude that $\{a_{j,k}^2\} \subset
L^{\infty}(\widetilde{\mathbb{H}}^+)$ as it is easy to verify
$x_1$ and $y_1$ are bounded with $N$.  Moreover, if we apply these
bounds to the relationship \eqref{4.1}, then we get the bounded
equation
\begin{equation}\label{4.2}
|g_j^{2}(z) - G_j^{2}(z)| \leq   K b_{1}^{(\omega(z) - \beta_2)} x_1 y_1, \quad \emph{for all} ~ z \in \widetilde{\mathbb{H}}^+,
\end{equation} where $K$ is an absolute constant which depends only upon the corona constants ($\mu$, $\delta,$ and $n$), the geometric considerations ($\epsilon$ and $\alpha_M)$, and our choice of $A$.  Lastly, under the map $\Phi^+(z)$ we can regard the newly constructed $\{g_j^2\}_{j=1}^n$ and $\{G_j^2\}_{j=1}^n$ as being functions defined on $\widetilde{\Omega}^+$.\\

\subsection*{The Subsequent Stitchings}

\indent \indent In the same fashion that we used to construct the relationship \eqref{4.1}, we  could construct  the third generation of solutions, $\{g_j^3\}_{j=1}^n$ and $\{G_j^3\}_{j=1}^n$, by stitching the newly formed $\{g_j^2\}_{j=1}^n$ to $\{g_j^1\}_{j=1}^n$ across the region $\mathcal{D}^-$ in $\widetilde{\mathbb{H}}^-$.  As soon as the third generation of solutions are constructed, we repeat the process, just as we did in the first stitchings, to obtain the fourth generation of solutions, $\{g_j^4\}_{j=1}^n$ and $\{G_j^4\}_{j=1}^n$, by stitching $\{g_j^3\}_{j=1}^n$ to $\{g_j^2\}_{j=1}^n$ across $\mathcal{D}^+$ in $\widetilde{\mathbb{H}}^+.$ Iterating this procedure with the sequences $\{b_m\}_{m=1}^{\infty}$, $\{x_m\}_{m=1}^{\infty}$, and $\{y_m\}_{m=1}^{\infty}$, we deduce the analogues of \eqref{4.1} and \eqref{4.2}:\\

\begin{align}
&g_j^m = G_j^m + \sum_{k=1}\left(a_{j,k}^m - a_{k,j}^m\right)f_k  &m=2,3, \ldots,  \\
&\left|g_j^{m+1}(z) - G_j^{m+1}(z)\right| \leq   K b_{m}^{(\omega(z) - \beta_2)} x_m y_m   &m=1,2, \ldots.  \label{4.4}
\end{align}\\

Since we will be referring to \eqref{4.4} many times from here, we
consider some variations. Each variation is customized to the
location of the variable $z$.  Recall,
\begin{equation*}
G_j^{m+1}(z) =
\begin{cases} g_j^m(z)     & \emph{when $m$ is odd, ~$~z \in \widetilde{\mathbb{H}}^+$, ~ and $z$ lies below $\delta_2^+$},
\\
g_j^{m-1}(z)       &\emph{when $m$ is odd, ~$~z \in \widetilde{\mathbb{H}}^+$, ~ and $z$ lies above $\delta_1^+$},
\end{cases}
\end{equation*}
\begin{equation*}
G_j^{m+1}(z) =
\begin{cases} g_j^m(z)     & \emph{when $m$ is even,~ $z \in \widetilde{\mathbb{H}}^-$, ~ and $z$ lies above $\delta_2^-$},
\\
g_j^{m-1}(z)       &\emph{when $m$ is even, ~$z \in \widetilde{\mathbb{H}}^-$, ~ and $z$ lies below $\delta_1^+$}.
\end{cases}
\end{equation*}\\

\noindent The phrases ``$z$ lies below $\delta_2^+$''  and ``$z$ lies above $\delta_1^+$'' when referring to $z \in \widetilde{\mathbb{H}}^+$ formally means $\omega(z)\geq \beta_2$ and $\omega(z) \leq \beta_1$ respectively.  Similarly when $z \in \widetilde{\mathbb{H}}^-$,  ``$z$ lies above $\delta_2^-$''  and ``$z$ lies below $\delta_1^+$''  means $\omega(z)\geq \beta_2$ and $\omega(z) \leq \beta_1$.  Immediately, we obtain two variations:
\begin{align*}
\left|g_j^{m+1}(z) - g_j^{m-1}(z)\right| &\leq   K b_{m}^{~- \beta_2} ~x_m y_m \quad \emph{ $m$ odd,~ $~z \in \widetilde{\mathbb{H}}^+$, and $z$ lies above $\delta_1^+$},\tag{4.4a} \label{4.4a} \\
\left|g_j^{m+1}(z) - g_j^{m-1}(z)\right| &\leq   K b_{m}^{~ - \beta_2}~ x_m y_m \quad \emph{ $m$ even, ~$z \in \widetilde{\mathbb{H}}^-$, and $z$ lies below $\delta_1^-$}. \tag{4.4b}
\end{align*}\\

We observe that in region \[ V^+ = \widetilde{\mathbb{H}}^+
\backslash \overline{\mathbb{H}^+} = \{z \in
\widetilde{\mathbb{H}}^+: ~ \emph{$z$ lies strictly below
$\bigcup_j I_j^+$ and strictly above
$\widetilde{\Psi}^+(\gamma_1^-)$} \}\]  we have a lower bound on
harmonic measure: $\omega(z) >1$.  With this observation and
\eqref{4.4}, we note
 \[ |g_j^{m+1}(z) - g_j^{m}(z)| \leq   K b_{m}^{(1- \beta_2)}
 ~x_m y_m \quad \emph{when $m$ is odd and $z \in V^+$}.\]
Under the map $\widetilde{\Psi}^- \circ \widetilde{\Phi}^+ $, we
transfer the preceding relationship to the extended lower half
plane (and repeat the construction for $m$ even with the region
$V^-$):

\begin{align*}
 |g_j^{m+1}(z) - g_j^{m}(z)| &\leq K b_{m}^{(1- \beta_2)} x_m y_m \quad \emph{when $m$ is odd, ~ $~z \in \widetilde{\Psi}^-(\widetilde{\Phi}^+ (V^+)) \subset \widetilde{\mathbb{H}}^-$}, \tag{4.4c} \\
 |g_j^{m+1}(z) - g_j^{m}(z)| &\leq K b_{m}^{(1- \beta_2)} x_m y_m \quad \emph{when $m$ is even,~ $z \in \widetilde{\Psi}^+(\widetilde{\Phi}^- (V^-)) \subset \widetilde{\mathbb{H}}^+$.} \tag{4.4d} \label{4.4d}
\end{align*}

\begin{wrapfigure}{l}{.5\textwidth}
\begin{center}
\includegraphics[width=.5\textwidth, bb= 99 74 691 535]{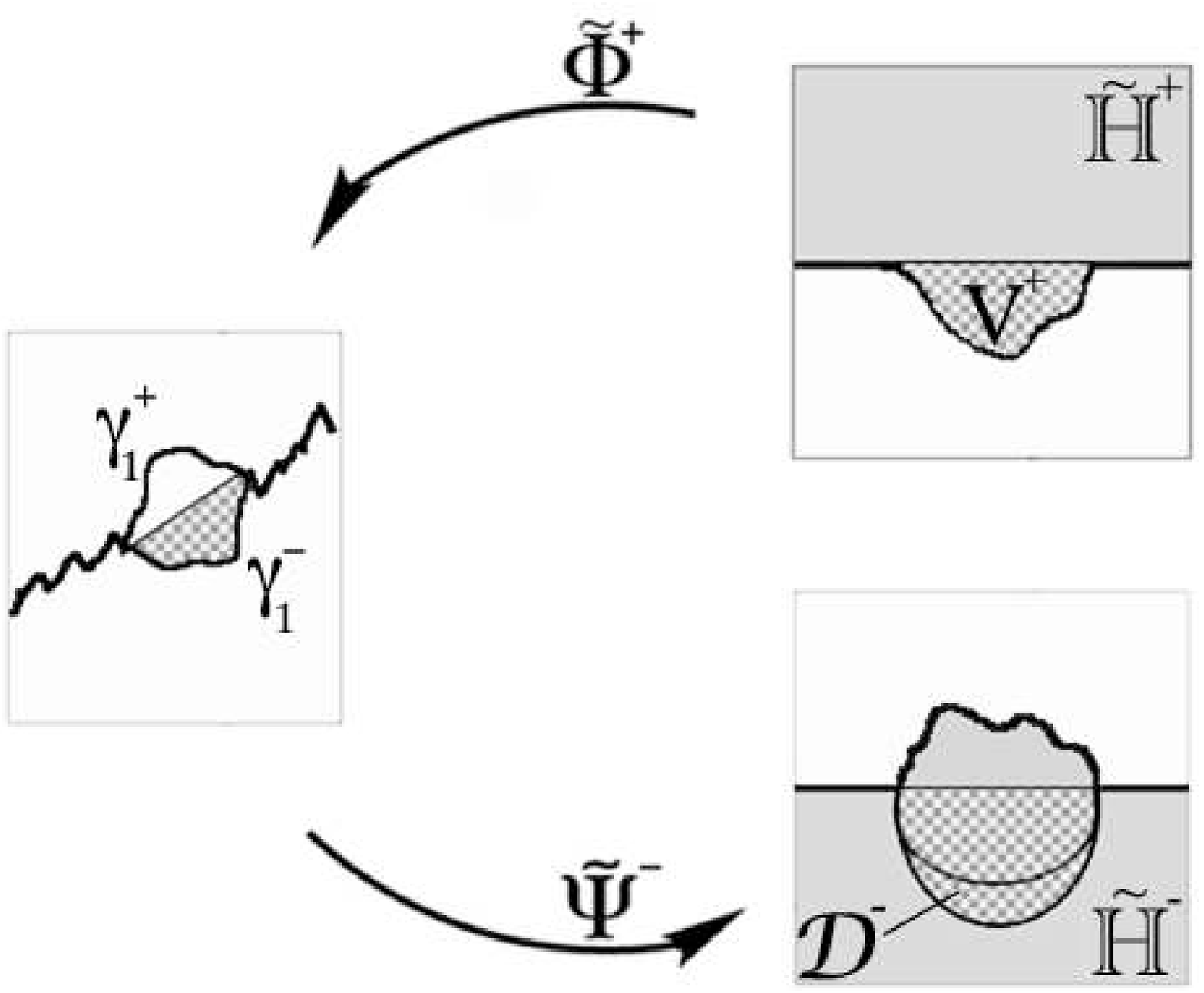}
\caption{\footnotesize The region $V+$ under the map $\widetilde{\Psi}^- \circ \widetilde{\Phi}^+.$}
\end{center}
\end{wrapfigure}

 From the latter  two variations and observing that $\mathcal{D}^+\subset \widetilde{\Psi}^+(\widetilde{\Phi}^- (V^-))$ and  $\mathcal{D}^-\subset \widetilde{\Psi}^-(\widetilde{\Phi}^+ (V^+))$, we deduce our first recursive relationship:\\

\noindent  \mbox{\qquad \quad  \: $x_1 \leq 2N$\quad}\\
\textbf{(R1)}  \\ \mbox{\qquad \quad  \: $x_{m+1} \leq K b_m^{~(1-\beta_2)} ~x_m y_m$}\\
\noindent \mbox{\qquad \qquad \qquad \qquad \qquad  \emph{for} $m=1,2,\ldots.$}\\

Next, let us deduce a recursive relationship for $\{y_m\}_{m=1}^{\infty}$.  It is easy to verify $y_1 \leq N$.  Now suppose $m$ is odd, then $y_{m+2} =\max_k\left\Vert G_k^{m+3} \right\Vert_{\mathcal{D}^+}$ $=\max_k \left\Vert g_k^{m+2}(\varphi^+) + g_k^{m+1}(1-\varphi^+)\right\Vert_{\mathcal{D}^+}.$\\

Let us take a look at the bounds for the functions in the above
equality.  Recall, \mbox{$\omega(z) > \beta_1$} when $z \in
\mathcal{D}^+$ so that \eqref{4.4} reduces to
\[ |g_k^{m+1}(z)| \leq |G_k^{m+1}(z)| + Kb_{m}^{~(\beta_1-\beta_2)}~ x_m y_m. \]
In addition, \eqref{4.4d} and the previous inequality imply
\begin{align*} |g_k^{m+2}(z)| &\leq |g_k^{m+1}(z)| + Kb_{m+1}^{(1-\beta_2)} x_{m+1} y_{m+1}\\
 &\leq |G_k^{m+1}(z)| + Kb_{m}^{(\beta_1-\beta_2)} x_m y_m + Kb_{m+1}^{(1-\beta_2)} x_{m+1} y_{m+1}.
 \end{align*}
Since the bound for $|g_k^{m+2}(z)|$ is greater than the bound for
$|g_k^{m+1}(z)|$, we deduce that
\begin{align*}
y_{m+2} &= \max_k \left\|g_k^{m+2}(\varphi^+) + g_k^{m+1}(1-\varphi^+)\right\|_{\mathcal{D}^+}\\
&\leq \max_k \left\||G_k^{m+1}(z)| + Kb_{m}^{(\beta_1-\beta_2)} x_m y_m + Kb_{m+1}^{(1-\beta_2)} x_{m+1} y_{m+1} \right\|_{\mathcal{D}^+}\\
&\leq \max_k \left\|G_k^{m+1}(z)\right\|_{\mathcal{D}^+} + Kb_{m}^{(\beta_1-\beta_2)} x_m y_m + Kb_{m+1}^{(1-\beta_2)} x_{m+1} y_{m+1}\\
&=  y_{m} + Kb_{m}^{(\beta_1-\beta_2)} x_m y_m + Kb_{m+1}^{(1-\beta_2)} x_{m+1} y_{m+1} .
\end{align*}

We could repeat verbatim  the case where $m$ is even, and for the
case of $y_2$ with the bound $\Vert g_j^1\Vert_{\infty} \leq N$,
so that we obtain the second recursive relationship:

\begin{align*}\label{R2}
&y_1 \leq N, \quad y_2\leq N+Kb_1^{(1-\beta_2)}~x_1y_1,\\
 \tag{R2}\\
  &y_{m+2} \leq y_{m} + Kb_{m}^{(\beta_1-\beta_2)} x_m y_m + Kb_{m+1}^{(1-\beta_2)} x_{m+1} y_{m+1} \qquad m=1,2, \ldots.
\end{align*}\\

\noindent The last two terms in this expression arise from the
error in stitching over two generations. It will be our goal to
show that these errors are summable.
\\
\\

\section{Convergence of \mbox{\boldmath${x_m}$}}
 In this chapter we select the factors
$\{b_m\}_{m=1}^{\infty}$ so that  $\{x_m\}_{m=1}^{\infty} \in \ell^1$ and $\{y_m\}_{m=1}^{\infty} \in \ell^{\infty}$ simultaneously:
\begin{lemma} \label{lemma5}
For positive constants $N$, $K$, $\beta_1$, and $\beta_2$, with
$\beta_1,\, \beta_2 <1$ and $(1-\beta_1)/(1-\beta_2) = 2$,
 there exists a sequence of positive real numbers $\{b_m\}_{m=1}^{\infty}$ with \[ 0<\inf_m\{b_m\} \qquad
 and \qquad b_m <1 \] such that for any
  pair of positive sequences $\{x_m\}_{m=1}^{\infty}$ and $\{y_m\}_{m=1}^{\infty}$ that satisfies the difference equations\\
\begin{align*}
    &x_{m+1} \leq K b_m^{~(1-\beta_2)} ~x_m y_m  \quad &m=1,2, \ldots&  \tag{R1} \label{R1}\\
  &y_{m+2} \leq y_{m} + Kb_{m}^{(\beta_1-\beta_2)} x_m y_m + Kb_{m+1}^{(1-\beta_2)} x_{m+1} y_{m+1} \quad &m=1,2, \ldots& \tag{R2}
  \end{align*}\\
  with initial data $x_1 \leq 2N,~ y_1 \leq N$, and $y_2 \leq N+Kb_1^{(1-\beta_2)}~x_1y_1$,
will also satisfy
 \[ \sum_{m=1}^{\infty} x_m < C_1<\infty \qquad and \qquad \sup_m\{y_m\}<C_2<\infty.\]\\
\end{lemma}

\noindent \textbf{Proof of Lemma 5.1:}\quad Since the pair of
sequences ($\{x_m\}_{m=1}^{\infty}$ and $\{y_m\}_{m=1}^{\infty}$)
that have equality holding in the initial data and equality
holding in \eqref{R1} and \eqref{R2} dominate all admissible
pairs, it suffices to solve for $\{b_m\}_{m=1}^{\infty}$ for this
particular pair. In addition, without loss of generality, we may
assume that $y_1, y_2 \geq 1$.  Now, fix $1>r>0$ (to be determined
later) and take $b_m$ so that
\begin{align*}
x_{m+1} & \stackrel{(R1)}{=}  K b_m^{(1-\beta_2)} x_m y_m = r^{m}  &\text{$m=1,2, \ldots.$} \quad \tag{5.1} \label{5.1}\\
\intertext{Next, we substitute the right hand side of the above equation to reduce \eqref{R2},}
 y_{m+2} &=   y_m + b_m^{(\beta_1 -1)} r^m    +r^{m+1}  &\text{$m=1,2, \ldots.$} \quad \tag{5.2} \label{5.2}
\end{align*}
\noindent When $m \geq 2$, we can solve for $b_m$ in terms of $r$
and $y_m$ by looking at successive generations of \eqref{5.1}.
Specifically, the right hand side of \eqref{5.1} at the $
m^{\small{th}}$ generation is
\begin{align*}
 &K b_m^{(1-\beta_2)} x_m y_m = r^m, \qquad \qquad \qquad \\
 \intertext{and by substituting the left hand side of \eqref{5.1} for $x_m$ makes}
&K b_m^{(1-\beta_2)} r^{m-1} y_m = r^m, \\
\text{so that} \qquad &K b_m^{\left(1-\beta_2\right)} y_m = r. \qquad \qquad \qquad
\end{align*}
If we raise both sides of the above equality to the power
$\left(\frac{\beta_1-1}{1-\beta_2}\right) = -2$, then
\begin{align*}
 &b_m^{\left(\beta_1-1\right)}  =  \left(\frac{Ky_m}{r}\right)^{2}. \tag{5.3a} \label{5.3a}\\
\intertext{For the case where $m=1$ we can repeat the preceding,
but with $x_1 = 2N$ to obtain}
&b_1^{\left(\beta_1-1\right)}  =  \left(\frac{Ky_1(2N)}{r}\right)^{2}. \tag{5.3b} \label{5.3b}
\end{align*}
If we substitute these relations for $b_m$ and $b_1$ into
\eqref{5.2}, then we have the ordinary difference equation:
\begin{align*}
y_1=N, \quad y_2=N+r, \quad y_{3} = y_1 + K^2 y_1^2 (2N)^2   r^{-1} + r^{2},  \tag{5.4a} \label{5.4a} \\
\notag \\
y_{m+2} = y_m + K^2~ y_m^2 ~r^{m-2} + r^{m+1} \qquad m= 2,3, \ldots  \tag{5.4b} \label{5.4b}
\end{align*}\\

To get bounds for $y_m$, we look at the difference
\begin{equation*}
 \frac{1}{y_m} - \frac{1}{y_{m+2}} = \frac{ y_{m+2} - y_{m}}{y_m y_{m+2}} = K^2 \left(\frac{y_m}{y_{m+2}} \right)r^{m-2}+ \frac{r^{m+1}}{y_m y_{m+2}}  \leq K^2r^{m-2} + r^{m+1}~ \quad m=2,3, \ldots
\end{equation*}
The last inequality holds since $1\leq \cdots \leq y_m \leq
y_{m+2}$.  By telescoping the differences, starting with $y_4$ and
$y_5$ for $m$ even and $m$ odd respectively,
\begin{align*}
 \frac{1}{y_4}-\frac{1}{y_{m+2}}  &\leq K^2 \sum_{\substack{j = 4,\\ j~ even}}^m r^{j-2}+ \sum_{\substack{j=4, \\ j~even}}^m r^{j+1} = \mathcal{O}(r^2)  &\emph{$m$ even, $m>2$},\\
 \frac{1}{y_5}-\frac{1}{y_{m+2}}  &\leq K^2 \sum_{\substack{j = 5,\\ ~j~ odd~}}^m r^{j-2}+ \sum_{\substack{j=5, \\ ~j~odd~}}^m r^{j+1} = \mathcal{O}(r^3)   &\emph{$m$ odd,  $~m>3$},
\intertext{We recall, $y_1 = N$, $y_2 = N+r$; and by using \eqref{5.4a} and
\eqref{5.4b}, $y_3 = \mathcal{O}(r^{-1})$, $y_4 = \mathcal{O}(1)$,
and $y_5 = \mathcal{O}(r^{-1})$.  So that for $r$ sufficiently small,}
0< C(r) < \frac{1}{y_4}& - \left(\ K^2 \sum_{\substack{j = 4,\\ j~ even}}^m r^{j-2}+ \sum_{\substack{j=4, \\ j~even}}^m r^{j+1}\right)\leq \frac{1}{y_{m+2}}  &\emph{$m$ even, $m>2$}, \\
\\
0< C(r) < \frac{1}{y_5}& - \left( K^2 \sum_{\substack{j = 5,\\ ~j~ odd~}}^m r^{j-2}+ \sum_{\substack{j=5, \\ ~j~odd~}}^m r^{j+1}\right)\leq \frac{1}{y_{m+2}}   &\emph{$m$ odd,  $~m>3$}.
\end{align*}
 Fix such an $r$ small enough so that the above inequalities holds and make sure
 $r>r_0>0$ so that,
\[\sup_m\{y_m\} <C_2(r_0) \qquad  \emph{and} \qquad  \displaystyle \sum_{m=1}^{\infty} x_m \leq 2N + \sum_{m=2}^{\infty} r^{m-1} <
C_1,\] while
\begin{align*}
&1> b_1 = \left(\frac{r}{Ky_1(2N)}\right)^{ \frac{2}{1-\beta_1}}\geq \left(\frac{r_0}{KC_2(2N)}\right)^{ \frac{2}{1-\beta_1}}\\
\intertext{and}
 &1> b_m =\left(\frac{r}{Ky_m}\right)^{\frac{2}{1-\beta_1}}\geq\left(\frac{r_0}{KC_2}\right)^{\frac{2}{1-\beta_1}}\quad m=2,3, \ldots
\end{align*}
We conclude that $\displaystyle \inf_m\{b_m\} > 0$, and thus
$\{b_m\}_{m=1}^{\infty}$ is our desired sequence.\qed
\\
\\

\section{Proof of Theorem 1.1.}

With the sequence $\{b_m\}_{m=1}^{\infty}$ following from Lemma
5.1, we now show:
\begin{align*}
\sup_{m~even} \left\Vert g_j^m\right\Vert_{H^{\infty}(\widetilde{\mathbb{H}}^+)} &\leq C<\infty \qquad j=1,\ldots, n, \\
 \emph{and} \quad \sup_{m~odd~} \left\Vert g_j^m\right\Vert_{H^{\infty}(\widetilde{\mathbb{H}}^-)} &\leq C<\infty \qquad j=1,\ldots, n,
\end{align*}
where $C$ is some absolute constant depending only upon
$\epsilon_0$, $\alpha$, and $A$.  We begin by looking at $\Vert
g_j^m\Vert_{H^{\infty}(\widetilde{\mathbb{H}}^+)}$ in the extended
upper half plane. Fix $z \in \widetilde{\mathbb{H}}^+$ and $m$
odd.  From \eqref{4.4} and variation \eqref{4.4a}  for the regions
$\{\omega(z)\leq \beta_1\}$, $\{\beta_1 < \omega(z) \leq
\beta_2\}$, and $\{\beta_2 < \omega(z)\}$ respectively, we have
\begin{align}
&|g_j^{m+1}(z) -g_j^{m-1}(z)| \leq Kb_{m}^{-\beta_2} x_{m} y_{m}   &\emph{when $z$ lies above $\delta_1^+$}&, \tag{6.1} \label{6.1}\\
&|g_j^{m+1}(z)| \leq |G_j^{m+1}(z)| + Kb_{m}^{(\beta_1 - \beta_2)}x_{m} y_{m} &\emph{when $z$ lies below $\delta_1^+$}& \emph{and above $\delta_2^+$,} \tag{6.2} \label{6.2} \\
&|g_j^{m+1}(z) - g_j^{m}(z)| \leq K x_{m} y_{m}  &\emph{when $z$ lies below $\delta_2^+$}&. \tag{6.3} \label{6.3}
\end{align}
Let us treat each region as its separate own special case.\\

\noindent \textbf{\textit{Case i})}: \textit{$z$ lies above $\delta_1^+$.}\\

 From the first relationship \eqref{6.1}, if we telescope
the differences over the even generations of $\{g_j^m\}$, then
\[ |g_j^{m+1}(z)| \leq N+ \sum_{\substack{k=1,\\k
~odd}}^{\infty} K b_{k}^{-\beta_2} x_k y_k \leq N+ K
\left(\frac{KC_2}{r_0}\right)^{
 \frac{2\beta_2}{1-\beta_1}} \sum_{\substack{k=1,\\k~odd}}^{\infty} x_k y_k \leq  C<\infty, \]\\
since $\sup_k\{y_k\} <C_2$, $\sum_k x_k < C_1$, and $\inf_k\{b_k\} \geq \left(\frac{r_0}{KC_2}\right)^{\frac{2}{1-\beta_1}}$ from Lemma 5.1.\\

\vspace{.3in}

\noindent \textbf{\textit{Case ii})} \textit{$z$ lies below $\delta_1^+$ and above $\delta_2^+$.}\\

For this case and the next we need some estimates similar to the
ones we obtained when we derived the second recursive relation,
\eqref{R2}.  Recall the relationship we have from \eqref{4.4d} in
this region,
\begin{align*}
&|g_j^1(z)-g_j^0(z)| \leq 2N, \\
\emph{and} \quad &|g_j^{m}(z) - g_j^{m-1}(z)| \leq K b_{m-1}^{(1-\beta_2)}~x_{m-1}y_{m-1} \qquad \emph{$m$ odd,} ~  m>1.
\intertext{As $G_j^{m+1}$ is an average of the two functions in the above,}
&|G_j^{m+1}(z)| = |g_j^{m}(z)(1-\varphi^+(z)) + g_j^{m-1}(z)\varphi^+(z)|,
\end{align*}
we can create two inequalities depending on whether we choose to
bound $g_j^{m+1}$ or $g_j^{m}:$
\begin{enumerate}
\item[\bf 1.)] $|G_j^{m+1}(z)| \leq
    |g_j^{m-1}(z)|+Kb_{m-1}^{(1-\beta_2)}x_{m-1}y_{m-1},$
    \quad \emph{m odd, ~ $m>1$}
\item[\bf 2.)] $|G_j^{m+1}(z)|\leq
    |g_j^{m}(z)|+Kb_{m-1}^{(1-\beta_2)}x_{m-1}y_{m-1}.$ \quad \emph{m odd, ~ $m>1$}
\end{enumerate}
If we choose the first inequality, \eqref{6.2} reduces to
\begin{align*}
&|g_j^2(z)| \leq 2N + Kb_1^{(\beta_1-\beta_2)}x_1 y_1, \tag{6.4a} \label{6.4a} \\
\emph{and} \quad &|g_j^{m+1}(z)| \leq |g_j^{m-1}(z)| + Kb_m^{(\beta_1 -\beta_2)}x_m y_m + Kb_{m-1}^{(1-\beta_2)}x_{m-1}y_{m-1}.
\intertext{and if we choose the second inequality \eqref{6.2} reduces to}
&|g_j^2(z)| \leq 2N + Kb_1^{(\beta_1-\beta_2)}x_1 y_1, \tag{6.4b} \label{6.4b}\\
\emph{and} \quad &|g_j^{m+1}(z)| \leq |g_j^{m}(z)| + Kb_m^{(\beta_1 -\beta_2)}x_m y_m + Kb_{m-1}^{(1-\beta_2)}x_{m-1}y_{m-1}.
\end{align*}
Now for Case \textit{ii}), \eqref{6.4a} unfolds to
\begin{align*}
 |g_j^{m+1}(z)| &\leq 2N + \sum_{\substack {k=1,\\ k~odd}}^{\infty} K b_k^{(\beta_1-\beta_2)}x_k y_k +\sum_{\substack {k=2, \\ k~even}}^{\infty} K b_{k}^{(1-\beta_2)}x_{k} y_{k}\\
\\
& \leq 2N + K\left(\frac{KC_2}{r_0}\right)^{ \frac{2(\beta_2 - \beta_1)}{1-\beta_1}} \sum_{k=1}^{\infty} x_k y_k \leq C <\infty.
\end{align*}

\vspace{.3in}

\noindent \textbf{\textit{Case iii})} \textit{$z$ lies below $\delta_2^+$.}\\

By the conformal map $\widetilde{\Psi}^+ \circ
\widetilde{\Phi}^-$, we have the relationships for the
$\mbox{\small{$m$-$\emph{1}$}}^{\mbox{\footnotesize{th}}}$ ($m$ odd, $m>1$)
generation
 of \eqref{6.3} and \eqref{6.4b} respectively:
\begin{align*}
&|g_j^{m}(z) - g_j^{m-1}(z)| \leq K x_{m-1} y_{m-1} \qquad
 \emph{when z lies above $\widetilde{\Psi}^+(\gamma_2^-)$},\\
 \\
&|g_j^{m}(z)| \leq |g_j^{m-1}(z)| + Kb_{m-1}^{(\beta_1 -\beta_2)}x_{m-1} y_{m-1} + Kb_{m-2}^{(1-\beta_2)}x_{m-2}y_{m-2}
 \quad \\
 &\emph{\hspace{2.4in}when z lies below $\widetilde{\Psi}^+(\gamma_2^-)$ and above $\widetilde{\Psi}^+(\gamma_1^-)$.}
\end{align*}\\
In the first event, we combine with \eqref{6.3} to get
\begin{align*}
|g_j^{m+1}(z)| &\leq |g_j^{m-1}(z)|+K x_{m}y_{m}+ Kx_{m-1}y_{m-1},\\
\intertext{and thus}
|g_j^{m+1}(z)| &\leq N+ \sum_{k=1}^{\infty} Kx_k y_k \leq C<\infty.
\end{align*}
 In the second event, we combine with \eqref{6.3} to get
\begin{align*}
 |g_j^{m+1}(z)| &\leq |g_j^{m-1}(z)|+ K x_{m} y_{m} +Kb_{m-1}^{(\beta_1 -\beta_2)}x_{m-1} y_{m-1} +Kb_{m-2}^{(1-\beta_2)}x_{m-2}y_{m-2},\\
 \intertext{and thus}
  |g_j^{m+1}(z)| &\leq N + 2K \sum_{k=1} b_k^{(\beta_1 - \beta_1)} x_ky_k \leq C <\infty.
  \end{align*}

 \vspace{.3in}

We conclude that $\Vert g_j^m\Vert_{H^{\infty}(\HH)} \leq C$ for all even $m$, since we have bounded these functions over the whole domain $\HH$.  Similarly, (with a lower bound) we could repeat the above procedure to conclude that $\Vert g_j^m\Vert_{H^{\infty}(\widetilde{\mathbb{H}}^-)} \leq C$ for all $m$ odd over the domain $\widetilde{\mathbb{H}}^-$.  Not only are $\{g_j^m\}_{\{m,even\}}$ and $\{g_j^m\}_{\{ m,odd\}}$ uniformly bounded in their respective domains, for each $j$, but also their difference has a shrinking bound in the intersection of $\widetilde{\Omega}^+$ and $\widetilde{\Omega}^-$.  We demonstrate this by showing that the sequences are uniformly Cauchy in $\Gamma \backslash E_0 = \bigcup_j F_j$ in the following sense:\\

Let $n>m\geq 0$ and let $z \in \Gamma \backslash E_0$, then with \eqref{4.4} and \eqref{5.1},
\begin{align}\label{6.5} \tag{6.5}
 |g_j^{n}(z) - g_j^{m}(z)| \leq \sum_{k =m}^{n-1} |g_j^{k}(z) - g_j^{k+1}(z)| \leq \sum_{k=m}^{n-1} Kb_k^{(1-\beta_2)}x_k y_k =\sum_{k=m}^{n-1} r^k. \end{align}\\

 Now, let $\{g_j^+\}_{j=1}^n$ be the normal limit of $\{g_j^m\}_{m=0}^{\infty} \subset H^{\infty}(\widetilde{\Omega}^+)$ for $m$ even, and let $\{g_j^-\}_{j=1}^n$ be the normal limit of $\{g_j^m\}_{m=1}^{\infty} \subset H^{\infty}(\widetilde{\Omega}^-)$ for $m$ odd.  As point-wise limits
  \[ g_1^{+}(z) f_1(z) + g_2^{+}(z) f_2(z) + \cdots + g_n^{+}(z) f_n(z) = 1 \qquad z \in \widetilde{\Omega}^+,\]
\[ g_1^{-}(z) f_1(z) + g_2^{-}(z) f_2(z) + \cdots + g_n^{-}(z) f_n(z) = 1 \qquad z \in \widetilde{\Omega}^-. \]

\noindent Moreover, \eqref{6.5} implies that $g_k^+(z) = g_k^-(z)$ for all $z \in \Gamma \backslash E_0$. Therefore, we can merge the two solutions together across $\Gamma \backslash E_0$, and obtain corona solutions on all of $\Omega$. \qed \\

For our proof, the homogeneous condition was critical.  Without it, we would not have been able to bound the crosscuts $\gamma_1^+$ and $\gamma_1^-$ into the disjoint diamonds, leaving the extended domains as multiply connected.  A proof for the non-homogeneous case still eludes the author. One might hope to avoid this obstacle by directly applying the results of the non-homogeneous cases, (e.g., the Denjoy domains).\\

The present work is part of the author's Ph.D. dissertation.  Most of all, the author would like to express his genuine gratitude to his thesis advisor, John Garnett, for countless hours of insightful conversations and guidance throughout the past couple years. The author is truly indebted for his support.


\newpage





\begin{thebibliography}{99}

\bibitem{alling1} Alling, N. L., A proof of the corona conjecture for finite open Riemann surfaces, \textit{Bull. Amer. Math. Soc.} \textbf{70} (1964), 110-112.\\

\bibitem{alling2} Alling, N. L., Extensions of meromorphic function rings over non-compact Riemann surfaces, I, \textit{Math. Z.} \textbf{89} (1965, 273-299.\\

\bibitem{behrens} Behrens, M., On the corona problem for a class of infinitely connected domains, \textit{Bull. Amer. Math. Soc.} \textbf{76} (1970), 387-391.\\

\bibitem{car-disc} Carleson, L., Interpolation by bounded analytic functions and the corona problem. \textit{Ann. of Math.} \textbf{76} (1962), 547-559.\\

\bibitem{car-fatsets} Carleson, L., On $H^\infty$ in multiply connected domains.  \textit{Conference on Harmonic Analysis in Honor of Antoni Zygmund}, Vol. 2, Wadsworth Inc., 1983, pp. 349-372.\\

\bibitem{earle-marden} Earle, C. J. and Marden, A., Projections to automorphic functions, \textit{Proc. Amer. Math. Soc.} \textbf{19} (1968), 274-278.\\

\bibitem{forelli} Forelli, F., Bounded holomorphic functions and projections, \textit{Illinois J. Math.} \textbf{10} (1966), 367-380.\\

\bibitem{gamelin} Gamelin, T.W., Wolff's proof of the Corona Theorem, \textit{Israel J. Math.} \textbf{37} (1980), 113-119.\\

\bibitem{gamelin2} Gamelin, T. W., Localization of the corona problem, \textit{Pac. J. Math.} \textbf{34} (1970), 73-81.\\

\bibitem{garnett} Garnett, J. B., \textit{Bounded Analytic Functions}, revised first edition, Graduate Texts in Mathematics 236, Springer, NY, 2007.\\

\bibitem{garnettjones} Garnett, J. B. and Jones, P. W., The Corona theorem for Denjoy domains, \textit{Acta. Math.} \textbf{155} (1985), 27-40.\\

\bibitem{gonzalez} Gonzalez, Maria Jose, Uniformly Perfect Sets, Green's function, and fundamental domains. \textit{Rev. Mat. Iberoamericana} \textbf{8} (1992), no. 2, 239-269.\\

\bibitem{gonzaleznicolau} González, María J. and Nicolau, Artur Quasiconformal mappings preserving interpolating sequences.  \textit{Ann. Acad. Sci. Fenn. Math.}  \textbf{23}  (1998),  no. 2, 283--290.\\

\bibitem{Handy} Handy, Jon, The Corona Theorem on Complements of
Certain Square Cantor Sets, \textit{arXiv 0712.1039} (2007). To
appear.\\

\bibitem{hormander} H\"ormander, L., Generators for some rings of analytic functions, \textit{Bull. Amer. Math. Soc.} \textbf{73} (1967), 943-949.\\

\bibitem{jones} Jones, P. W., Carleson measures and the Fefferman-Stein
decomposition of BMO (R). \textit{Ann. of Math.}, \textbf{111}
(1980), 197-208.\\

\bibitem{jones2}  ----- $L^\infty$  estimates for the $\overline{\partial}$ problem in a half-plane. \textit{Acta
Math.}, \textbf{150} (1983), 137-152.\\

\bibitem{jones3} ---- \textit{Some problems in complex analysis,}  The Bieberbach Conjecture, Proceedings of the Symposium on the Occasion of the Proof The Bieberbach conjecture (West Lafayette, Ind., 1985) Math. Surverys Monographs,  21, Amer. Math. Soc., Providence, R.I. (1986), 105-108.\\

\bibitem{jones-marshall} Jones, P. W. and Marshall, D. E., Critical points of Green's function, harmonic measure, and the corona problem, \textit{Ark. f\"or Mat.} \textbf{23} (1985), no. 2, 281-314.\\

\bibitem{Kenig} Kenig, Carlos E., Weighted $H^p$ spaces on
Lipschitz Domains. \textit{American Journal of Mathematics}, Vol.
102, No. 1 (Feb., 1980), 129-163.\\

\bibitem{Muckenhoupt} Muckenhoupt, B., The equivalence of two conditions for weight functions. \textit{Studia Math.} \textbf{49} (1974), 101-106.\\

\bibitem{Pommerenke} Pommerenke, Ch., Uniformly Perfect Sets and
the Poincar\'{e} metric. \textit{Arch. Math.} \textbf{32} (1979),
192-199.\\

\bibitem{slod2} Slodkowski, Z., On bounded analytic functions in finitely connected domains, \textit{Tran. Amer. Math. Soc.} \textbf{300} (1987), 721-736.\\

\bibitem{stout1} Stout, E. L., Two theorems concerning functions holomorphic on multiply connected domains, \textit{Bull. Amer. Math. Soc.} \textbf{69} (1963), 527-530.\\

\bibitem{stout2} Stout, E. L., Bounded holormorphic functions on finite Riemann surfaces, \textit{Trans. Amer. Math. Soc.} \textbf{120} (1965), 255-285.\\

\bibitem{stout3} Stout, E. L., On some algebras of analytic functions on finite open Riemann surfaces, \textit{Math. Z.} \textbf{92} (1966), 366-379.  Corrections to: On some algebras of analytic functions on finite open Riemann surfaces, \textit{Math. Z.} \textbf{95} (1967), 403-404.\\

\bibitem{Tsuji} Tsuji,M., \textit{Potential Theory In Modern Function Theory}. Maruzen Co., LTD. TOKYO.1959.\\

\bibitem{varopo} Varopoulos, N.  Th., Ensembles pics et ensembles d'interpolation d'une algebra uniforme, \textit{C. R. Acad. Sci. Paris, Ser. A.} \textbf{272} (1971), 866-867.\\


\end{thebibliography}
\end{document}